\documentclass{amsart}

\usepackage{amsfonts}
\usepackage{amsmath}
\usepackage{amsthm}
\usepackage{amssymb}
\usepackage[english]{babel}
\usepackage{graphics}
\usepackage{latexsym}
\usepackage{longtable} 
\usepackage{mathrsfs} 
\usepackage{graphicx}
\usepackage{wrapfig} 
\usepackage[pdftex,colorlinks=true,linkcolor=blue,citecolor=green]{hyperref}
\usepackage{enumitem}
\usepackage{esint}
\usepackage{mathcomp}
\usepackage{stmaryrd}
\usepackage{wasysym}
\usepackage{fancyhdr}

\newtheorem{theorem}{Theorem}
\newtheorem{corollary}[theorem]{Corollary}
\newtheorem{lemma}[theorem]{Lemma}
\newtheorem{proposition}[theorem]{Proposition}

\newtheorem{claim}[theorem]{Claim}
\newtheorem{example}[theorem]{Example}

\theoremstyle{definition}
\newtheorem{definition}[theorem]{Definition}

\renewcommand{\d}{\mathrm{d}}

\newcommand{\mL}{\mathcal{L}}
\newcommand{\mH}{\mathcal{H}}

\newcommand{\mX}{\mathfrak{X}}

\newcommand{\mM}{\mathcal{M}}

\newcommand{\E}{\textbf{E}}
\newcommand{\M}{\mathrm{M}}
\newcommand{\F}{\mathrm{F}}
\newcommand{\D}{\mathrm{D}}

\newcommand{\X}{\mathfrak{X}}

\newcommand{\A}{\textbf{A}}

\newcommand{\K}{\mathrm{K}}
\newcommand{\G}{\mathrm{G}}

\newcommand{\R}{\mathbb{R}}

\newcommand{\N}{\mathbb{N}}

\newcommand{\mB}{\mathbb{B}}

\newcommand{\Q}{\mathrm{Q}}

\newcommand{\noi}{\noindent}
\newcommand{\ms}{\medskip}

\newcommand{\al}{\alpha}

\newcommand{\ga}{\gamma}
\newcommand{\Ga}{\Gamma}

\newcommand{\e}{\varepsilon}
\newcommand{\si}{\sigma}

\newcommand{\la}{\lambda}
\newcommand{\La}{\Lambda}

\newcommand{\Om}{\Omega}
\newcommand{\om}{\omega}

\newcommand{\av}{-\hspace{-10.5pt}\displaystyle\int}

\newcommand{\weak }{\, -\!\!\!\!\!-\!\!\!\!\rightharpoonup}
\newcommand{\weakstar }{ \overset{\, *_{\phantom{|}}}{{\smash{\, -\!\!\!\!-\!\!\!\!\rightharpoonup}}\, } }

\newcommand{\larrow}{\longrightarrow}
\newcommand{\ot}{\otimes}

\newcommand{\LL}{\text{\LARGE$\llcorner$}}
\newcommand{\p}{\partial}
\newcommand{\sub}{\subseteq}
\newcommand{\set}{\setminus}
\newcommand{\by}{\times}

\newcommand{\ess}{\mathrm{ess}}

\newcommand{\dist}{\mathrm{dist}}

\renewcommand{\div}{\mathrm{div}}

\newcommand{\bt}{\begin{theorem}}\newcommand{\et}{\end{theorem}}
\newcommand{\bd}{\begin{definition}}\newcommand{\ed}{\end{definition}}
\newcommand{\bl}{\begin{lemma}}\newcommand{\el}{\end{lemma}}
\newcommand{\beq}{\begin{equation}}\newcommand{\eeq}{\end{equation}}
\newcommand{\bc}{\begin{claim}}\newcommand{\ec}{\end{claim}}
\newcommand{\bex}{\begin{example}}\newcommand{\eex}{\end{example}}
\newcommand{\bcor}{\begin{corollary}}\newcommand{\ecor}{\end{corollary}}
\newcommand{\bp}{\begin{proof}}\newcommand{\ep}{\end{proof}}

\newcommand{\BPP}{\medskip \noindent \textbf{Proof of Proposition} }
\newcommand{\BPT}{\medskip \noindent \textbf{Proof of Theorem} }

\numberwithin{equation}{section}

\makeatletter
\@namedef{subjclassname@2020}{%
  \textup{2020} Mathematics Subject Classification}
\makeatother

\begin{document}

\title[ISOSUPREMIC VECTORIAL MINIMISATION PROBLEMS IN $L^\infty$ WITH  CONSTRAINTS]{On isosupremic vectorial minimisation problems in $L^\infty$ with general nonlinear constraints}

\author{Ed Clark}

\address{E. C., Department of Mathematics and Statistics, University of Reading, Whiteknights Campus, Pepper Lane, Reading RG6 6AX, United Kingdom}

\email{e.d.clark@pgr.reading.ac.uk}
 
\author{Nikos Katzourakis}

\address[Corresponding author]{N. K., Department of Mathematics and Statistics, University of Reading, Whiteknights Campus, Pepper Lane, Reading RG6 6AX, United Kingdom}

\email{n.katzourakis@reading.ac.uk}

\thanks{E.C.\ has been financially supported through the UK EPSRC scholarship GS19-055}
  
\subjclass[2020]{35J94; 49J20; 49K20; 49K35; 49K99; 35J92.}

\date{}

\keywords{Calculus of Variations in $L^\infty$; PDE-Constrained Optimisation; Euler-Lagrange equations; Absolute minimisers; Aronsson-Euler systems; Kuhn-Tucker theory; Isosupremic problems, Isoperimetric problems.}

\begin{abstract} We study minimisation problems in $L^\infty$ for general quasiconvex first order functionals, where the class of admissible mappings is constrained by the sublevel sets of another supremal functional and by the zero set of a nonlinear operator. Examples of admissible operators include those expressing pointwise, unilateral, integral isoperimetric, elliptic quasilinear differential, jacobian and null Lagrangian constraints. Via the method of $L^p$ approximations as $p\to \infty$, we illustrate the existence of a special $L^\infty$ minimiser which solves a divergence PDE system involving certain auxiliary measures as coefficients. This system can be seen as a divergence form counterpart of the Aronsson PDE system which is associated with the constrained $L^\infty$ variational problem.
\end{abstract}

\maketitle

\tableofcontents

\section{Introduction and main results}   \label{Section1}

Let $n, N \in \N$ and let also $\Om \Subset \R^n$ be a bounded open set with Lipschitz continuous boundary. Consider two functions
\beq 
\label{1.1}
  f,\, g\ : \  \ \Om\times \R^N \times \R^{N\times n} \longrightarrow \R,
\eeq
which will be assumed to satisfy certain natural structural assumptions. Additionally, let $\bar p>n$ be fixed and consider a given nonlinear operator
\beq \label{1.2}
\begin{split}
\Q\ :\ \ W^{1, \bar p}_0\big(\Om;\R^N\big) \longrightarrow \E,
\end{split}
\eeq
where $(\E,\|\cdot\|)$ is an arbitrary Banach space. In this paper we are interested in the following variational problem: given $G\geq 0$ and the supremal functionals
\[
\mathrm{F}_{\infty}, \mathrm{G}_{\infty}\ :\ \ W^{1,\infty}_0(\Om;\R^N) \larrow \R,
\]
defined by
\beq
\label{1.3}
\left\{ \ \ 
\begin{split}
\mathrm{F}_{\infty}(u) :&= \, \underset{\Om}{\ess \sup} \, f(\cdot, u, \D u) ,\phantom{\Big|}
\\
\mathrm{G}_{\infty}(u) :&=\, \underset{\Om}{\ess \sup} \, g(\cdot, u, \D u) , 
\end{split}
\right. \ \
\eeq
find $u_\infty \in W^{1,\infty}_0(\Om;\R^N)$ such that
\beq
\label{1.4}
\mathrm{F}_{\infty}(u_\infty)\,=\, \inf \Big\{ \mathrm{F}_{\infty}(u)\ : \  u\in W^{1,\infty}_0(\Om;\R^N),\ \mathrm{G}_{\infty}(u)\leq G\ \,\& \, \  \Q(u)=0\Big\}.
\\
\eeq
We are also interested in deriving appropriate differential equations as necessary conditions that such constrained minimisers $u_\infty$ might satisfy. Hence, the problem under consideration is that of minimising the supremal functional $\textrm{F}_\infty$ within the $G$-sublevel set of the supremal functional $\textrm{G}_\infty$, constrained also by the zero level set of the mapping $\Q$. Minimisation problems in $L^\infty$ with an $L^\infty$ sublevel-set constraints are called \emph{isosupremic}, a terminology introduced by Aronsson-Barron in \cite{AB}, by analogy to the terminology of \emph{isoperimetric} constraints utilised for integral functionals.

The assumptions required for our mapping $\Q$ expressing the nonlinear level-set constraint are very general, allowing for the solvability of a very large class of problems. In Section \ref{section5} we provide a fairly comprehensive list of explicit examples of operators $\Q$ to which our results apply. Indeed, $\Q$ may manifest itself as a pointwise, unilateral, or inclusion constraint (Subsection \ref{S4.1}), as an integral isoperimetric constraint (Subsection \ref{S4.2}), or as a nonlinear PDE constraint, including second order quasilinear divergence systems (Subsection \ref{S4.3}), jacobian equations and other PDEs driven by null Lagrangians (Subsection \ref{S4.4}).  

Minimisation of supremal functionals  poses several important challenges. 
Firstly, the $L^\infty$ norm is neither Gateaux differentiable nor strictly convex, and the space $L^\infty$ is neither reflexive nor separable. A related additional complication is that the $L^\infty$ norm is not additive but only sub-additive with respect to the domain. As a result, if one solves the $L^\infty$ minimisation existence problem via the direct method, then the analogue of the Euler-Lagrange equations for the $L^\infty$ problem cannot be derived directly by considering variations, due to the lack of smoothness of the $L^\infty$ norm.

In this paper, we transcend the difficulties illustrated above, by following the relatively standard strategy of using appropriate approximations by $L^p$ functionals as $p\to \infty$. We solve the desired $L^\infty$ variational problem by solving approximating $L^p$ variational problems for all $p$, and obtaining the necessary compactness estimates that allow us to pass to the limit as $p\to \infty$. The case of finite $p$, which is of independent interest and a byproduct of our analysis, for us is just the mechanism to solve the desired $L^\infty$ problem. The intuition behind the use of $L^p$ approximations is based on the fact that for a fixed essentially bounded function on a set of finite measure, the $L^p$ norm tends to the $L^\infty$ norm as $p\to\infty$. 

Problems involving constraints are relatively new in the Calculus of Variations in $L^\infty$ and previous work has been relatively sparse, even more so in the vectorial case. To the best of our knowledge, the only work which directly studies isosupremic problems is \cite{AB} by Aronsson-Barron. Among other questions answered therein, it considers some aspects of the one-dimensional case for $n=1$, but with no additional constraints of any type (which amounts to $\Q\equiv 0$ in our setting). 

More broadly, very few previous works involve vectorial problems with general constraints in $L^\infty$. Certain vectorial and higher order problems involving eigenvalues in $L^\infty$ have been considered in \cite{K4, KP}. Examples of problems with PDE and other constraints are considered in \cite{CKM, K1, K2, K3}. In the paper \cite{BJ} of Barron-Jensen, a scalar $L^\infty$ constrained problem was considered, but the constraint was integral. With the exception of the paper \cite{AB}, it appears that vectorial variational problems in $L^\infty$ involving isosupremic constraints have not been studied before, especially including additional nonlinear constraints which cover numerous different cases, as in this work. For assorted interesting works within the wider area of Calculus of Variations in $L^\infty$ we refer to \cite{AP, BBJ, BJW1, BJW2, BN, BP, CDP, KM, KZ, MWZ, PZ, RZ}.

Let us note that, in this work, we refrain from discussing the question of defining and studying localised versions of $L^\infty$ minimisers on subdomains. Such minimisers are commonly called \emph{absolute} in the Calculus of Variations in $L^\infty$, a concept first introduced and studied by Aronsson in his seminal papers in the 1960s. As already noted in \cite{AB}, the definition of constrained absolute minimisers on subdomains is quite problematic in many respects, even for specific choices of constraints in the scalar case $N=1$, which has a well developed $L^\infty$ theory. In any case, the ``non-intrinsic" method of $L^p$ approximations, employed herein, is widely believed to select always the ``best" $L^\infty$ minimiser, which in the scalar unconstrained case can be shown to be indeed an absolute minimiser. In the vectorial case of $N\geq 2$, the situation is trickier and the ``correct" localised minimality notion is still under discussion, even without any constraints being involved (see e.g.\ \cite{AK1, AK2} for work in this direction).
 
We now state our main results. Our notation is generally standard as e.g.\ in \cite{KV}, or otherwise self-explanatory. As we will only be working with finite measures, we will be using the normalised $L^p$ norms in which we replace the integral over the domain with the respective average. For our density functions $f$ and $g$ appearing in \eqref{1.1} we will assume the following:
\beq 
\label{1.5}
\left\{
\begin{split}
&\text{Both $f$ and $g$ are Carath\'eodory functions, namely $f(\cdot,\eta,P)$, $g(\cdot,\eta,P)$ are}
\\
& \text{measurable for all $(\eta,P)$, and $f(x,\cdot,\cdot)$, $g(x,\cdot,\cdot)$ are continuous for a.e.\ $x$.}
\end{split}
\right. \!\!\!\!
\eeq
\beq
\label{1.6}
\left\{
\begin{split}
& \text{Exists a continuous $\mathrm{C} : \overline{\Om}\times \R^N \larrow [0,\infty)$ and an $\alpha>1$:}
\\
\ms
& \hspace{40pt} \left\{ \ \ 
\begin{split}
0\leq f(x,\eta,P)\leq \mathrm{C}(x,\eta)\big(|P|^{\alpha}+1\big),& 
\\
0\leq g(x,\eta, P) \leq \mathrm{C}(x,\eta)\big(|P|^{\alpha}+1\big), 
\end{split}
\right.
\\
& \text{for a.e.\ $x$ and all $(\eta,P)$.}
\end{split}
\right.
\eeq 
\beq 
\label{1.7}
\text{For a.e.\ $x$ and all $\eta$, $f(x,\eta,\cdot)$ and $g(x,\eta,\cdot)$ are quasiconvex on $\R^{N\by n}$.}
\eeq
\beq
\label{1.8}
\left\{
\begin{split}
& \text{\texttt{Either} $f$ \texttt{or} $g$ is coercive, namely exist $c,C>0$ such that either}
\\
&\hspace{80pt} f(x, \eta, P) \, \geq \, c|P|^{\alpha}- C, 
\\
 &\text{or}  
 \\
&\hspace{80pt} g(x, \eta, P) \, \geq \, c|P|^{\alpha}- C,
\\
& \text{for a.e.\ $x$ and all $(\eta,P)$.}
\end{split}
\right.
\eeq
These assumptions are relatively standard in the Calculus of Variations. We note only that \eqref{1.7} is meant in the sense of Morrey quasiconvexity for integral functionals (as e.g.\ in \cite{D}), not in the sense of level-convexity or of the $L^\infty$ notion of ``BJW-quasiconvexity" of Barron-Jensen-Wang in \cite{BJW2}, nor in the sense of ``$\mathcal{A}$-Young quasiconvexity" of Ansini-Prinari in \cite{AP}. This stronger notion of quasiconvexity is not necessary if one is interested in merely solving \eqref{1.4} by applying the direct method in $W^{1,\infty}$ \emph{without deriving any PDEs}, and can indeed be weakened substantially, but it is needed for Theorems \ref{theorem2} and \ref{theorem3}, so we simplify the exposition by assuming it at the outset. For our operator $\Q$ in \eqref{1.2} we will assume the following:
\beq 
\label{1.9}
\Q^{-1}\big(\{0\}\big) \ \text{is weakly closed  in} \ W^{1,\overline {p}}_0(\Om;\R^N).
\eeq
This is a very feeble requirement for $\Q$ and there exist numerous explicit examples of interest that satisfy \eqref{1.9}, see Section \ref{section5}. For $1\leq p<\infty$, we define the approximating $L^p$ functionals $\mathrm{F}_p , \mathrm{G}_p : W^{1,\al p}_0(\Om;\R^N) \larrow \R$ by setting
\beq
\label{1.10}
\left\{ 
\ \ \ \
\begin{split}
\mathrm{F}_{p}(u) :&=  \left(\, {\av}_{\!\!\!\Om} f(\cdot,u,\D u)^p \, \mathrm{d} \mL^{n}\! \right)^{\!\!1/p},\\
\mathrm{G}_{p}(u) :&=  \left(\, {\av}_{\!\!\!\Om} g(\cdot,u,\D u)^p \, \mathrm{d} \mL^{n}\! \right)^{\!\!1/p}.
\end{split}
\right.
\eeq 
For $p=\infty$, $\mathrm{F}_\infty$ and $\mathrm{G}_\infty$ are given by \eqref{1.3}. For each $p \in [1,\infty]$ and $G\geq 0$ fixed, we define the admissible minimisation class $\mathfrak{X}^{p}({\Om})$ on which $\mathrm{F}_p$ is to be minimised, by setting
\beq
\label{1.11}
\mathfrak{X}^p(\Om) \,:=\, \Big\{ v \in W^{1,\al p}_0(\Om;\R^N) \ : \ \mathrm{G}_p(v)\leq G \ \text{and} \ \Q(v)=0\Big\}.
\eeq

Our first main result concerns the existence of $\mathrm{F}_p$-minimisers in $\X^p(\Om)$ and the existence of $\mathrm{F}_\infty$-minimisers in $\X^\infty(\Om)$, obtained as subsequential limits as $p\to \infty$. 

\begin{theorem}[$\mathrm{F}_\infty$-minimisers, $\mathrm{F}_p$-minimisers \& convergence as $p\to\infty$] \label{theorem1} 
Suppose that the mappings $f, g$ and $\Q$ satisfy the assumptions \eqref{1.5} through \eqref{1.9}. If the next compatibility condition is satisfied
\beq 
\label{1.13}
\inf \Big\{ \G_{\infty} : \Q^{-1}\big(\{0\}\big) \cap W^{1,\infty}_0(\Om;\R^N) \Big\} \,<\, G,
\eeq
then, for any $p\in[\bar p,\infty]$, the functional $\mathrm{F}_p$ has a constrained minimiser $u_p$ in the admissible class $\mathfrak{X}^p(\Om)$, namely
\beq \label{1.14}
\mathrm{F}_p\big(u_p)\, =\, \inf\Big\{\mathrm{F}_p\big(v)\, : \ v \in \mX^p(\Om) \Big\}.
\eeq
Additionally, there exists a subsequence of indices $(p_j)_1^\infty$ such that, the sequence of respective $\mathrm{F}_{p_j}$-minimisers satisfies $u_{p} \larrow u_\infty$ uniformly on $\overline{\Om}$, and $u_{p} \weak u_\infty$ weakly in $W^{1, q}_0(\Om;\R^N)$, for all $q\in (1,\infty)$ fixed, as $p_j\to\infty$. Finally, we have the convergence of minimum values $\mathrm{F}_p (u_{p}) \larrow  \mathrm{F}_\infty (u_\infty)$ as $p_j\to\infty$.
\end{theorem}

Note that, in view of \eqref{1.9}, the compatibility condition \eqref{1.13}  guarantees that $\mX^p(\Om)\neq \emptyset$. Given the existence of constrained minimisers provided from Theorem \ref{theorem1} above, the next natural question concerns the deduction of necessary conditions satisfied by constrained minimisers in the form of PDEs. Firstly, we examine the case of $p<\infty$. The nonlinear constraint expressed by the zero level-set of the mapping $\Q$ in \eqref{1.11} gives rise to a Lagrange multiplier in the Euler-Lagrange equations. This can be inferred by employing well-known results on Lagrange multipliers, see for instance \cite{Z}. For this to be possible, though, one needs improved regularity of the mappings $f, g$ and $\Q$ involved. We will suppose additionally that
\beq
\label{1.15}
\left\{ \ \ 
\begin{split}
& \text{The partial derivatives $f_{\eta}, f_P, g_{\eta} , g_P$ of $f$ and $g$ are continuous on \ \ \ }
\\
&\text{$\smash{\overline{\Om}\times \R^N \times \R^{N\times n}}$, and for $\mathrm{C}, \al $ as in \eqref{1.6}, we have the bounds}
\\
& \ \ \ \ \ \left\{ \ \
\begin{split}
& |f_\eta|(x,\eta, P)+|f_P|(x,\eta, P) \leq \mathrm{C}(x,\eta)\big(|P|^{\alpha-1}+1\big),
\\
&   |g_\eta|(x,\eta, P)  +|g_P|(x,\eta, P) \leq \mathrm{C}(x,\eta)\big(|P|^{\alpha-1}+1\big),
\end{split}
\right.
\\
&\text{for all }(x,\eta,P).
\end{split}
\right.
\eeq
Further, we will assume that:
\beq
\label{1.16}
\left\{
\begin{split}
&\text{$\Q$ is continuously differentiable, and its Fr\'echet derivative}
 \\
&\hspace{50pt} (\d \Q)_{\overline{u}}\ :\ \  W^{1, \bar p}_0(\Om;\R^N)\longrightarrow \E
\\
&\text{has closed range in $\E$, for any $\overline{u} \in \Q^{-1}(\{0\})\sub W^{1, \bar p}_0(\Om;\R^N)$. }
\end{split}
\right.
\eeq
Recall that no regularity was assumed for $\Q$ to obtain the existence of minimisers in Theorem \ref{theorem1}. Finally, for the sake of brevity, for any $u \in W^{1,1}_{\text{loc}}(\Om)$, we will employ the following notation
\[
f[u]\,\equiv \, f(\cdot, u, \D u), \ \ \ g[u]\,\equiv \, g(\cdot, u, \D u),
\]
and similar notation will be used for the compositions of $f_{\eta}, f_P, g_{\eta} , g_P$ of $f$ and $g$ respectively, namely $f_\eta[u], g_\eta[u]$, etc. Further, ``$\cdot$" and ``$:$" will denote the standard inner products on $\R^N$ and $\R^{N\by n}$ respectively. 

\begin{theorem}[The equations in $L^p$] \label{theorem2} 
Suppose we are in the setting of Theorem \ref{theorem1} and assumptions \eqref{1.5} through \eqref{1.9} are satisfied. Suppose additionally that \eqref{1.15} and \eqref{1.16} are satisfied. Then, for any $p\in(\bar p,\infty)$, there exists Lagrange multipliers
\beq
\label{1.17}
\lambda_p \geq 0 , \ \ \mu_p\geq 0 , \ \ \ \psi_p \, \in \, \emph{\E}^*,
\eeq
where $(\emph \E^*,\|\cdot\|_*)$ is the dual space of $ \emph\E$, such that not all vanish simultaneously:
\beq
\label{1.18}
|\lambda_p| + |\mu_p| + \|\psi_p\|_* \ne 0.
\eeq
Then, the minimiser $u_{p} \in \mX^p(\Om)$ satisfies the equation
\beq
\label{1.19}
 \left\{ 
\begin{split}
 & \lambda_p \, {\av}_{\!\!\!\Om} f[u_p]^{p-1}\Big(f_{\eta}[u_p] \cdot \phi \, +\, f_P[u_p] : \D \phi\Big) \, \d \mL^n 
 \\
+ \ & \mu_p \, {\av}_{\!\!\!\Om} g[u_p]^{p-1} \Big(g_{\eta}[u_p]\cdot \phi \, +\, g_P[u_p] : \D \phi \Big) \, \d {\mathcal{L}}^{n} \, =\, \big\langle \psi_p , (\d \Q)_{u_p}(\phi) \big\rangle,
\end{split}
\right.
\eeq
for all test mappings $\phi \in W^{1, \al p}_0\big(\Om;\R^N\big)$, coupled by the additional condition
\beq
\label{1.18}
\mu_p \big(\G_p(u_p) - G\big)=0.
\eeq
\end{theorem}

Note that condition \eqref{1.18} implies that if $\G_{p}(u_p)<G$, namely if the $L^\infty$-energy constraint is not realised (i.e.\ the minimiser $u_p$ lies in the interior of the sublevel set $\{\G_p \leq G\}$), then $\mu_p=0$ and hence the associated multiplier vanishes.

Now we consider the case of $p=\infty$. In this case, there is not a simple analogue of the divergence structure Euler Lagrange equations.  The equations are derived by an appropriate passage to limits as $p\to\infty$ in Theorem \ref{theorem2}. A standard approach in Calculus of Variations in $L^\infty$ has been to derive Aronsson-type PDE systems, which are non-divergence counterparts to the Euler-Lagrange equations, as e.g.\ done in \cite{AB} for the case of $n=1$ and $\Q\equiv 0$. However, Aronsson-type systems are always non-divergence and far less tractable than their divergence counterparts. In fact, in the vectorial case they have discontinuous coefficients and are fully nonlinear in the higher order case (see for instance \cite{CKP} and \cite{K-1, K0, KPr} for this evolving line of development regarding the direct study of generalised solutions to Aronsson systems).

Nevertheless, there exists an alternative approach which allows to derive \emph{divergence structure} PDE systems as necessary conditions. The starting point of this idea is based on the use of a different scaling in the Euler-Lagrange equations in $L^p$ and has already born substantial fruit in \cite{CKM, K1, K2, K3, K4, KM, KP}. There is however a toll to be paid for this ``forcing" of divergence structure: certain non-uniquely determined measures arise as auxiliary parameters in the coefficients of the PDE system, which depend nonlinearly on the minimisers. For more details on the historical origins of this alternative approach to deriving $L^\infty$ equations for variational problems we refer to \cite{K4}.

Our final main result therefore concerns the satisfaction of necessary conditions for the constrained minimiser in $L^\infty$ constructed in Theorem \ref{theorem1}. For this result we will need to impose some natural additional hypotheses. These hypotheses, although they restrict considerably the classes of $f, g, \Q$ that were utilised in order to prove existence of minimisers, they do nonetheless include the interesting case of $F_\infty$ being the $L^\infty$ norm of the gradient. Firstly, let us introduce some convenient notation and rewrite \eqref{1.19} in a way which will be more appropriate for the statement and the subsequent proof. By introducing for each $p \in (\bar p,\infty)$ the non-negative Radon measures $\si_p,\tau_p \in \mM(\overline{\Om})$ given by
\beq
\label{1.20}
\begin{split}
\si_p\, &:=\, \left\{
\begin{array}{ll}
\dfrac{1}{\mL^n(\Om)}\bigg(\dfrac{f[u_p]}{ \F_p(u_p)}\bigg)^{\! p-1} \mL^{n}\LL_{\Om}, \ \ & \text{ if }\F_p(u_p)>0,
\\
0, \ \ & \text{ if }\F_p(u_p)=0,
\end{array}
\right. 
\\
\tau_p\, &:=\, \left\{
\begin{array}{ll}
\dfrac{1}{\mL^n(\Om)}\bigg(\dfrac{g[u_p]}{ \G_p(u_p)}\bigg)^{\! p-1} \mL^{n}\LL_{\Om}, \ \ & \text{ if }\G_p(u_p)>0,
\\
0, \ \ & \text{ if }\G_p(u_p)=0,
\end{array}
\right.
\end{split}
\eeq
and the scaled multipliers
\beq
\label{1.21}
\begin{split}
\hat\la_p\, &:=\, \left\{
\begin{array}{ll}
\la_p \F_p(u_p)^{p-1}, \ \ & \text{ if }\F_p(u_p)>0,
\\
\la_p, \ \ & \text{ if }\F_p(u_p)=0,
\end{array}
\right. 
\\
\hat \mu_p\, &:=\, \left\{
\begin{array}{ll}
\mu_p \G_p(u_p)^{p-1}, \ \ & \text{ if }\G_p(u_p)>0,
\\
\mu_p, \ \ & \text{ if }\G_p(u_p)=0,
\end{array}
\right.
\end{split}
\eeq
we can rewrite \eqref{1.19} as
\beq
\label{1.22}
 \left\{ 
\begin{split}
 & \hat\lambda_p \int_\Om \Big(f_{\eta}[u_p] \cdot \phi \, +\, f_P[u_p] : \D \phi\Big) \, \d \si_p 
 \\
+ \ & \hat \mu_p \int_\Om \Big(g_{\eta}[u_p] \cdot \phi \, +\, g_P[u_p] : \D \phi \Big) \, \d \tau_p \, =\, \big\langle \psi_p , (\d \Q)_{u_p}(\phi) \big\rangle,
\end{split}
\right.
\eeq
Further, let us set
\beq
\label{1.23}
R_p \,:=\, \hat\la_p + \hat\mu_p + \|\psi_p\|_*
\eeq
and note that, by virtue of Theorem \ref{theorem2} and Definition \eqref{1.21}, we have that $R_p>0$. We may then define the new rescaled multipliers
\beq
\label{1.24}
\La_p \,:=\, \frac{\hat\lambda_p}{R_p} \in[0,1], \ \ \ \M_p \,:=\, \frac{\hat\mu_p}{R_p} \in[0,1], \ \ \ \Psi_p \,:=\, \frac{\psi_p}{R_p} \in \bar\mB^{\E^*}_1\!(0),
\eeq
where $\bar \mB^{\E^*}_1\!(0)$ is the closed unit ball in $\E^*$. Let us finally set
\[
C^1_0\big(\overline{\Om};\R^N \big) \,:=\, W^{1,\infty}_0(\Om;\R^N) \bigcap C^1\big(\R^n;\R^N \big).
\]
Now we state the additional assumptions which we need to impose:
\beq
\label{1.25}
\E \text{ is a separable Banach space}.
\eeq
\beq
\label{1.26}
\left\{ \ \ 
\begin{split}
&\text{The restriction of the differential $(u,v)\mapsto (\d \Q)_{u}(v)$, considered as}
\\
& \hspace{60pt} \d \Q \ :\ \  \Q^{-1}(\{0\}) \by W^{1, \bar p}_0(\Om; \R^N) \longrightarrow \E, 
\\
&\text{satisfies the following conditions:}
\\
& \text{If }u_m \weak u \text{ in }\Q^{-1}(\{0\}) \text{ as }m\to\infty, \text{ and }\phi \in W^{1, \bar p}_0(\Om; \R^N), \text{ then }
\\
& \hspace{60pt} 
\left\{
\ \ 
\begin{split}
 & (\d \Q)_{u_m}(u_m)  \larrow (\d \Q)_u (u), 
\\
& (\d \Q)_{u_m}(\phi) \larrow (\d \Q)_u (\phi),
\end{split}
\right.
\\
&\text{ as }m\to\infty. 
\end{split}
\right.
\eeq
The above assumption requires that $\d \Q$ be weakly-strongly continuous on the diagonal of $\Q^{-1}(\{0\}) \by \Q^{-1}(\{0\})$ and on subsets of the form $\Q^{-1}(\{0\}) \by \{\phi\}$, when $W^{1, \bar p}_0(\Om; \R^N) \by W^{1, \bar p}_0(\Om; \R^N)$ is endowed with its weak topology and $\E$ with its norm topology. We assume further that: 
\beq
\label{1.27}
\left\{ \ \
\begin{split}
&\text{(i) $g$ does not depend on $P$, namely $g(x,\eta,P)=g(x,\eta)$,}
\\
&\text{(ii) $f$ is quadratic in $P$ and independent of $\eta$, namely}
\\
&\text{\hspace{80pt} $f(x,\eta,P)=\A(x):P\ot P$,}
\\
&\text{for some continuous positive symmetric fourth order tensor} 
\\
&\text{$\smash{\A : \overline{\Om}\larrow \R^{N\by n} \ot \R^{N\by n}}$, which satisfies} 
\\
&\text{\hspace{30pt}$\A(x):P\ot P>0$\ , \ \ $\A(x):P\ot Q = \A(x) : Q\ot P$,}
\\
&\text{for all $x \in \overline{\Om}$ and all $P,Q \in \R^{N\by n}\set\{0\}$.}
\end{split}
\right.
\eeq
The above requirements are compatible with the previous assumptions on $f$. In fact, by \cite[Lemma 4, p. 8]{K4} and our earlier assumptions, the positivity and symmetry requirements for $\A$ are superfluous and can be deduced by merely assuming that $f$ is quadratic in $P$ (up to a replacement of $\A$ by its symmetrisation), but we have added them to \eqref{1.27} for simplicity. We may finally state our last principal result. 

\begin{theorem}[The equations in $L^\infty$] \label{theorem3} 
Suppose we are in the setting of Theorem \ref{theorem2} and that the same assumptions are satisfied. Additionally we assume that \eqref{1.25} through \eqref{1.27} hold true. Then, there exist
\beq
\label{1.28}
\La_\infty \in [0,1], \ \ \ \emph M_\infty \in [0,1], \ \  \Psi_\infty \in \bar \mB^{\emph \E^*}_1\!(0),
\eeq
which are Lagrange multipliers associated with the constrained minimisation problem \eqref{1.4}. There also exist Radon measures
\beq
\label{1.29}
\begin{split}
\sigma_\infty \in \mM ({\overline{\Om}} ), \ \ \ \tau_\infty \in \mM (\overline{\Om} ),
\end{split}
\eeq
and a Borel measurable mapping $\D u_\infty^\star : \overline{\Om}\larrow \R^{N\by n}$ which is a version of $\D u_\infty \in L^\infty(\Om;\R^{N\by n})$, such that the minimiser $u_\infty \in \mX^\infty(\Om)$ satisfies the equation
\beq
\label{1.30}
\begin{split}
 \La_\infty \int_{\overline{\Om}}  f_P(\cdot, \D u_\infty^\star) : \D \phi  \, \d \si_\infty \, +\, \emph M_{\infty} \int_{\overline{\Om}} g_{\eta}(\cdot, u_\infty) \cdot \phi \, \d \tau_\infty
  \, =\, \big\langle \Psi_\infty , (\d \Q)_{u_\infty}(\phi) \big\rangle,
\end{split}
\eeq
for all test maps $\phi\in C^1_0\big(\overline{\Om};\R^N \big)$, coupled by the condition
\beq
\label{1.31}
\emph M_\infty \big(\G_\infty(u_\infty) - G\big)=0. 
\eeq
Additionally, the map $\D u_\infty^\star$ can be represented as follows:
\beq
\label{1.32}
\left\{ \ \ 
\begin{array}{l}
\text{For any sequence $(v_j)_1^\infty \sub C^1_0\big(\overline{\Om};\R^N \big)$ satisfying} 
\ms\ms
\\
\ \ \ \ \ \ \ \ \ \ \left\{  
\begin{split}
& \lim_{j\to\infty} \|v_j - u_\infty\|_{(W^{1,1}_0 \cap L^\infty)(\Om)} =\, 0,
\\
 & \ \ \ \limsup_{j\to\infty} \F_\infty(v_j)  \leq \,  \F_\infty(u_\infty),   
\end{split}
\right.
\ms
\\
\text{exist a subsequence $(j_k)_1^\infty$ such that} 
\ms\smallskip
 \\ 
 \D u_\infty^\star (x) \, = \left\{
\begin{array}{ll}
\underset{k \to \infty}{\lim}\, \D v_{j_k} (x), & \text{if the limit exists},
\\
0, & \text{otherwise}.
\end{array}
\right.
\smallskip
\end{array}
\right.
\eeq
(Such an explicit sequence $(v_j)_1^\infty$ is constructed in the proof.) Finally, the Lagrange multipliers $\La_{\infty}, \emph M_\infty, \Psi_\infty$ and the measures $\sigma_\infty, \tau_{\infty}$ can be approximated as follows:
\beq
\label{1.33}
\left\{ \ \
\begin{split}
\Psi_p \, &\weakstar  \, \Psi_\infty,  \ \ \text{ in } \bar \mB^{\emph \E^*}_1\!(0),
\\
\La_p \,  &\larrow  \, \La_\infty, \ \,  \text{ in } [0,1], 
 \\
\mathrm M_p \,  &\larrow  \, \mathrm M_\infty,  \   \text{ in } [0,1], 
\end{split}
\right.
\eeq 
and
 \beq
\label{1.34}
\left\{ \ \
\begin{split}
\si_p \, \weakstar  \, \si_\infty, \ \  &\text{ in }\mM (\overline{\Om}  ),
\\
\tau_p  \, \weakstar \,  \tau_\infty,  \ \ \  &\text{ in }\mM (\overline{\Om} ),\end{split}
\right.
\eeq
along a subsequence $p_j \to \infty$.
\end{theorem}

The weak interpretation of \eqref{1.30} is
\[
- \La_\infty  \div(f_P(\cdot, \D u_\infty^\star) \si_\infty\big)  \, +\, \mathrm M_\infty  \, g_{\eta}(\cdot,u_\infty) \tau_\infty  \, =\, \big\langle \Psi_\infty , (\d \Q)_{u_\infty} \big\rangle,
\]
in $\big( C^1_0\big(\overline{\Om};\R^N \big) \big)^*$, up to the identifications
\[
\big\langle \Psi_\infty , (\d \Q)_{u_\infty} \big\rangle \equiv \big\langle \Psi_\infty , (\d \Q)_{u_\infty} (\cdot) \big\rangle, \ \ \ g_\eta \equiv g_\eta \cdot (\cdot), \ \ \ f_P \equiv (\cdot) \cdot f_P.
\]
Note that in Theorem \ref{theorem3}, the equations obtained depend on certain measures not a priori known explicitly. Therefore, their significance is understood to be largely theoretical, rather than computational. For the proof of this result, we will utilise some machinery developed in the recent paper \cite{K4} for some related work on generalised $\infty$-eigenvalue problems. The main points of this approach are recalled in the course of the proof, for the convenience of the reader.

We conclude this lengthy introduction with some comments concerning the composition of this paper. In Sections \ref{section2} and \ref{section3} we establish our main results, Theorems \ref{theorem1}, \ref{theorem2} and \ref{theorem3}. In Section \ref{section5} we provide a rather detailed list of explicit large classes of nonlinear operators $\Q$ to which our results apply.

\section{Minimisers of $L^p$ problems and convergence as $p\to\infty$}  \label{section2}

In this section we demonstrate Theorem \ref{theorem1}. The proof is a consequence of the next two propositions, utilising the direct method of the Calculus of Variations.

\begin{proposition}\label{Proposition 4}
In the setting of Theorem \ref{theorem1} and under the same assumptions, for any $p \in [\bar p, \infty)$, the functional $\mathrm{F}_p$ has a constrained minimiser $u_p \in \mathfrak{X}^p(\Om)$, as claimed in \eqref{1.14}.
\end{proposition}

\BPP \ref{Proposition 4}.  Fix $p \geq \bar p >n$. We begin by illustrating that $\mathfrak{X}^p(\Om)\ne \emptyset $. First note that by the compatibility condition \eqref{1.13}, the finiteness of the infimum implies that $\Q^{-1}\big(\{0\}\big) \cap W^{1,\infty}_0(\Om;\R^N) \neq \emptyset$. Further, there exists $u_0\in W^{1,\infty}_0(\Om;\R^N) $ with $\Q(u_0)=0$ such that $\mathrm{G}_{\infty}(u_0) <G$. Hence, by H\"older inequality we have
\[
\begin{split}
\mathrm{G}_{p}(u_0) \,= \, \left(\, {\av}_{\!\!\!\Om}g(\cdot, u_0, \D u_0)^p\, \mathrm {d} \mL^{n}\! \right)^{\!\!1/p} \leq \, \big\|g(\cdot, u_0, \D u_0) \big\|_{L^{\infty}(\Om)} \, =\, \mathrm{G}_{\infty}(u_0) \, < \, G.
\end{split}
\]
Consequently, in view of \eqref{1.11}, both constraints are satisfied by $u_0$, hence $u_0\in \mathfrak{X}^p(\Om) \neq \emptyset$. Next, note that $f^p$ is a (Morrey) quasiconvex function. To see this, let $h:\R^{N\times n}\larrow \R$ be an arbitrary quasiconvex function, in our case we will take $h(P)=f(x, \eta, P)$ for fixed $(x,\eta)$. Then, by assumption \eqref{1.7}, for any $\phi\in W^{1,\infty}_0(U; \R^N)$ with $U\Subset \R^n$ open and $P\in \R^{N \by n}$ fixed,
$$
h(P)\leq \av_{U} h(P+\D\phi) \ \d \mathcal{L}^n.
$$
Hence, by Jensen's inequality and the convexity of $t\mapsto t^p$, we conclude that
$$
h(P)^p\leq \bigg(\, {\av}_{\!\!\!U} h(P+\D\phi) \, \d \mathcal{L}^n \!\bigg)^{\!\!p} \, \leq \, {\av}_{\!\!\! U} h(P+\D\phi)^p \, \d \mathcal{L}^n.
$$
We now proceed to bound $f^p$. By \eqref{1.6}, we estimate
\[
\begin{split}
&0 \leq f(x, \eta, P)^p\leq \mathrm{C}(x, \eta)^p \big( 1+ |P|^{\alpha}\big)^p\leq 2^{p-1} \mathrm{C}(x, \eta)^p \big( 1+ |P|^{\alpha p}\big).
\end{split}
\]
By standard results (see \cite{D}), $\mathrm{F}_p$ is weakly lower semicontinuous on $W^{1,\alpha p}_0(\Om;\R^N)$. Let $(u^{(i)})_1^{\infty} \sub \mathfrak{X}^p(\Om)$ denote a minimising sequence. By virtue of  \eqref{1.6} we have $f \geq 0$, hence clearly $\inf_{i\in\N} \mathrm{F}_p(u^{(i)}) \geq 0$. We now show that the infimum is finite. To this aim, by \eqref{1.6} we estimate
\[
\begin{split}
\inf _{i\in\N} \mathrm{F}_p(u^{(i)}) &\leq \mathrm{F}_p(u_0)
\\
&=\bigg(\, {\av}_{\!\!\!\Om} f(\cdot,u_0, \D u_0)^p \ \d{\mathcal{L}}^n \bigg)^{\!\! 1/p}
\\
&\leq \bigg(\, {\av}_{\!\!\!\Om} \Big(\mathrm{C}(\cdot, u_0)(1+|\D u_0|^{\alpha})\Big)^p \, \d {\mathcal{L}}^n \bigg)^{\!\!1/p},
\end{split}
\]
which yields
\[
\begin{split}
\inf _{i\in\N} \mathrm{F}_p(u^{(i)}) & 
\leq \Big\|\mathrm{C}(\cdot, u_0)\big(1+|\D u_0|^{\alpha}\big) \Big\|_{L^{\infty}(\Om)}
\\
&\leq \|\mathrm{C}(\cdot, u_0)\|_{L^{\infty}(\Om)}\|(1+|\D u_0|^{\alpha})\|_{L^{\infty}(\Om)}
\\
&\leq \|\mathrm{C}(\cdot, u_0)\|_{L^{\infty}(\Om)}\big(1+ \| \D u_0\|^\al_{L^{\infty}(\Om)}\big)
\\
&<\infty.
\end{split}
\]
Hence, the infimum is indeed finite. Now we show that under assumption \eqref{1.8}, the minimising sequence $(u^{(i)})_1^{\infty}$ is bounded in $W^{1,\al p}(\Om;\R^N)$. Let $h$ symbolise either $f$ or $g$, whichever is coercive. Hence, since $h$ is coercive, we have
\[
{\av}_{\!\!\!\Om} \big|h(\cdot, u^{(i)}, \D u^{(i)}) + C \big|^p \, \d \mathcal{L}^n \, \geq \, c \ {\av}_{\!\!\!\Om} |\D u^{(i)}|^{\alpha p} \ \d \mathcal{L}^n.
 \]
By using the Poincar\'e and H\"older inequalities, we infer that
\[
\begin{split}
&C\, +\, \bigg(\, {\av}_{\!\!\!\Om} |h(\cdot, u^{(i)}, \D u^{(i)})|^p  \ \d \mathcal{L}^n\bigg)^{\!\!\frac{1}{p}} \, \geq \, c'\, \|u^{(i)}\|^{\alpha}_{W^{1,\alpha p}(\Om)}
\end{split}
\]
for some new constant $c'>0$ which is independent of $i \in\N$. If $h=f$, then by the previously derived estimates we have the uniform bound
\[
 \|u^{(i)}\|^{\alpha}_{W^{1,\alpha p}(\Om)} \, \leq\, \frac{1}{c'}\Big(C \,+\, \F_p(u^{(i)}) \Big) \, \leq\, \frac{C}{c'} \,+\, \frac{1}{c'}\|\mathrm{C}(\cdot, u_0)\|_{L^{\infty}(\Om)}\big(1+ \| \D u_0\|^\al_{L^{\infty}(\Om)}\big), 
\]
and if $h=g$, then by the isosupremic constraint we have the uniform bound
\[
 \|u^{(i)}\|^{\alpha}_{W^{1,\alpha p}(\Om)} \, \leq\, \frac{1}{c'}\Big(C \,+\, \G_p(u^{(i)}) \Big) \, \leq\, \frac{C \,+\, G}{c'} .
\]
In either case, we have that $(u^{(i)})_1^{\infty}$ is weakly precompact in $W^{1,\al p}_0(\Om;\R^N)$. By passing to a subsequence if necessary, standard strong and weak compactness arguments imply that there exists a map $u_p \in W^{1,\al p}_0(\Om;\R^N)$ and a subsequence denoted again by $(u^{(i)})_1^{\infty}$ such that
\[
\left\{ \ 
\begin{split}
 &u^{(i)} \larrow u_p,  \ \ \ \ \ \ \,\text{in} \ L^{\alpha p}(\Om; \R^N), \smallskip\\
 & \D u^{(i)} \weak \D u_p,  \ \ \text{in} \ L^{\alpha p}(\Om; \R^{N\times n}), \smallskip\\
\end{split}
\right.
\]
as $i\to \infty$. Further, since $p>n$, by the Morrey estimate we have that $(u^{(i)})_1^{\infty}$ is also bounded in $C^{0,\gamma}(\overline{\Om}, \R^N)$ for $\ga<1-n/(\al p)$. By the compact embedding of H\"older spaces, we conclude that
\[
u^{(i)}\larrow u_p \ \ \ \text{in}\ C(\overline{\Om}; \R^N),
\]
as $i\to \infty$. It remains to show that $\mathfrak{X}^p(\Om)$ is weakly closed in $W^{1,\al p}_0(\Om;\R^N)$. To this end, we need to show that the constraints determine a weakly closed subset of $W^{1, \alpha p}_0(\Om; \R^N)$. Firstly note that by assumptions \eqref{1.5}-\eqref{1.7}, $\mathrm{G}_p$ is a weakly lower semi-continuous functional on $W^{1, \alpha p}_0(\Om; \R^N)$. This can be seen by an analogous argument to that used to show that $\mathrm{F}_p$ is weakly lower-continuous. Since $(u^{(i)})_1^{\infty} \sub \mathfrak{X}^p(\Om)$, we have $\mathrm{G}_p(u^{(i)})\leq G$ \ for all $i\in \N$. Therefore, since $u^{(i)} \weak u_p$ in $W^{1,\al p}_0(\Om;\R^N)$ as $i\to\infty$ and $\mathrm{G}_p$ is weakly lower-continuous, we infer that
\[
\begin{split}
\mathrm{G}_p(u_p) \, \leq \, \liminf_{i\to\infty} \mathrm{G}_p(u^{(i)})\, \leq \, G.
\end{split}
\]
Hence, $\mathrm{G}_p(u_p)\leq G$. Recall now that, in view of assumption \eqref{1.9}, $\Q^{-1}(\{0\})$ is a weakly closed subset of $W^{1,\bar p}(\Om;\R^N)$, where $\bar p>n$. We now show that $\Q^{-1}(\{0\})$ is also a weakly closed subset of $W^{1,\al p}_0(\Om;\R^N)$, which will complete the claim that $\mathfrak{X}^p(\Om)$ is weakly closed in $W^{1,\al p}_0(\Om;\R^N)$. Note first that since $(u^{(i)})_1^{\infty} \sub \mathfrak{X}^p(\Om)$, we have that $\Q(u^{(i)})=0$ for all $i\in\N$. Recalling that $u^{(i)} \weak u_p$ in $W^{1,\al p}_0(\Om;\R^N)$ as $i\to\infty$ and that $\al p > \bar p$ because by assumption $\al>1$, we deduce that $u^{(i)} \weak u_p$ in $W^{1,\bar p}(\Om;\R^N)$ as well, as $i\to\infty$. Since $\Q^{-1}(\{0\})$ is weakly closed in $W^{1,\bar p}(\Om;\R^N)$, we infer that $\Q(u_p)=0$, as desired. Thus, for each $p\in [\bar p,\infty)$, $u_p \in \mathfrak{X}^p(\Om)$, and
\[
\F_p(u_p) \, \leq \,\liminf_{i\to\infty}\F_p(u^{(i)})\, =\, \inf\big\{\F_p\, :\ \mathfrak{X}^p(\Om)\big\}. 
\]
The proposition ensues.
\qed

\ms

Our next result below establishes the existence of constrained minimisers for $\mathrm{F}_{\infty}$ and the approximation by minimisers of the $\mathrm{F}_p$ functionals as $p\to \infty$, therefore completing the proof of Theorem \ref{theorem1}.

\begin{proposition}\label{Proposition 5}
In the setting of Theorem \ref{theorem1} and under the same assumptions, the functional $\mathrm{F}_{\infty}$ has a constrained minimiser $u_\infty \in \mathfrak{X}^\infty(\Om)$, as claimed in \eqref{1.14}. Additionally, the claimed modes of convergence  hold true for a subsequence of minimisers $(u_{p_j})_{j=1}^\infty$ as $j\to\infty$.
\end{proposition}

\BPP \ref{Proposition 5}. We continue from the proof of Proposition \ref{Proposition 4}. The element $u_0\in\mathfrak{X}^p(\Om)$ provides an energy bound uniform in $p$, and also $u_0\in \mathfrak{X}^{\infty}(\Om)$. Fix $p,q>1$ with $p \geq q\geq\bar p$. By the H\"older inequality, minimality and the definition of the constrained class, we have the estimates
\[
\left\{ \ \ 
\begin{split}
\mathrm{F}_{q}(u_p) \, & \leq \, \mathrm{F}_p(u_p) \, \leq \, \mathrm{F}_p(u_0)
\, \leq \, \|\mathrm{C}(\cdot, u_0)\|_{L^{\infty}(\Om)}\big(1+\|\D u_0\|^{\alpha}_{L^{\infty}(\Om)}\big),
\\
\mathrm{G}_{q}(u_p) \, &\leq \, \mathrm{G}_p(u_p)\leq G,
\end{split}
\right.
\]
with right hand side bounds which are uniform in $p, q$. We now argue in a similar fashion to that used in the proof of Proposition \ref{Proposition 4}. We first show that under assumption \eqref{1.8}, the family of minimisers $(u_p)_{p\geq \bar p}$ is bounded in $W^{1,q}(\Om;\R^N)$, for any $q\in (1,\infty)$ fixed. Let $h$ symbolise either $f$ or $g$, whichever is coercive. We then we have that
\[
{\av}_{\!\!\!\Om} \big|h(\cdot, u_p, \D u_p) + C \big|^q \, \d \mathcal{L}^n \, \geq \, c \ {\av}_{\!\!\!\Om} |\D u_p|^{\alpha q} \ \d \mathcal{L}^n.
 \]
Since $\al>1$, by the H\"older inequality we infer
\[
C\, +\, \bigg(\, {\av}_{\!\!\!\Om} |h(\cdot, u_p, \D u_p)|^q  \ \d \mathcal{L}^n\bigg)^{\!\!\frac{1}{q}} \, \geq \, c^{\frac{1}{q}}\, \|\D u_p\|^{\alpha}_{L^{q}(\Om)},
\]
for the constants $c,C>0$ of \eqref{1.8} (which are independent of $p$ and $q$). If $h=f$, by applying our earlier estimates we deduce the uniform bound
\[
 \|\D u_p\|^{\alpha}_{L^{q}(\Om)} \, \leq\, \frac{1}{c}\Big(C \,+\, \F_q(u_p)\Big)  \, \leq\, \frac{C}{c} \,+\, \frac{1}{c}\|\mathrm{C}(\cdot, u_0)\|_{L^{\infty}(\Om)}\big(1+ \| \D u_0\|^\al_{L^{\infty}(\Om)}\big).
\]
If $h=g$, then again as in our earlier estimates we have the uniform bound
\[
 \|\D u_p\|^{\alpha}_{L^{q}(\Om)} \, \leq\, \frac{1}{c}\Big(C \,+\, \G_q(u_p)\Big)  \, \leq\, \frac{C \,+\, G}{c}.
\]
In either case, we see that under \eqref{1.8}, our estimates above imply that 
\[
 \|\D u_p\|_{L^{q}(\Om)} \, \leq\, K,
\]
for some constant $K>0$ independent of $p,q$. Further, by the Poincar\'e inequality, we deduce that
\[
 \|u_p\|_{W^{1,q}(\Om)} \, \leq\, K\big(1+C(q)\big),
\]
where $C(q)$ is the constant of the Poincar\'e inequality in $W^{1,q}_0(\Om;\R^N)$.
Hence, the sequence of minimisers $(u_p)_{p\geq \bar p}$ is bounded in $W^{1,q}_0(\Om;\R^N)$ for any fixed $q\in (1,\infty)$, and therefore it is weakly precompact in this scale of spaces. By a standard diagonal argument, there exists a sequence $(p_j)_1^\infty$ and a mapping  
$$ 
u_\infty \, \in \, \bigcap_{\bar p <q<\infty} W^{1, q}_0(\Om; \R^n),  
$$ 
such that $u_p \weak u_\infty$ in $W^{1, q}_0(\Om; \R^n)$ as $p_j \to \infty$, for any fixed $q\in (\bar p, \infty)$. By standard compactness arguments in Sobolev and H\"older spaces, we infer that
\[
\left\{ \ 
\begin{split}
& u_p\larrow u_{\infty}, \ \ \ \ \ \ \, \text{in}\ C\big(\overline{\Om}; \R^N\big),\\
 & \D u_p \weak \D u_{\infty},  \ \ \text{in} \ L^{q}(\Om; \R^{N\times n}), \smallskip\\
\end{split}
\right.
\]
as $p_j\to\infty$, for any $q\in (\bar p, \infty)$. We will now show that $u_\infty\in \mathfrak{X}^{\infty}(\Om)$. In view of \eqref{1.11}, we need to show that $u_\infty \in W^{1,\infty}_0(\Om; \R^N)$ and that $\G_\infty(u_\infty)\leq G$ and also $\Q(u_\infty)=0$. By the weak lower semi-continuity of the $L^q$ norm for $q\geq \bar p$ fixed, we have
\[
\|\D u_\infty\|_{L^{q}(\Om)} \, \leq\, \liminf_{p_j \to \infty} \|\D u_p\|_{L^{q}(\Om)} \, \leq\, K.
\]
By letting $q\to \infty$, this yields that $\D u_\infty \in L^{\infty}(\Om; \R^N)$. By the Poincar\'e inequality in $W^{1,\infty}_0(\Om;\R^N)$, we infer that $u_\infty \in W^{1,\infty}_0(\Om; \R^N)$. Next, since $\mathrm{G}_p(u_p)\leq G$ for all $p\in(\bar p, \infty)$, via the H\"older inequality and weak lower semi-continuity, we have
\[
\begin{split}
\mathrm{G}_{\infty}(u_{\infty})= \lim_{q\to \infty}\mathrm{G}_{q}(u_{\infty})
&\leq \liminf_{q\to \infty} \Big(\liminf_{p_j\to \infty}\mathrm{G}_{q}(u_{p})\Big)
\, \leq \, \liminf_{p_j\to \infty}\mathrm{G}_{p}(u_{p})
\, \leq \, G,
\end{split}
\]
yielding that indeed $\mathrm{G}_{\infty}(u_{\infty}) \leq G$. We now show that $\Q(u_\infty)=0$. We have already shown in Proposition \ref{Proposition 4} that $\Q^{-1}(\{0\})$ is a weakly closed subset of $W^{1,q}_0(\Om;\R^N)$ for any $q\in (\bar p,\infty)$. Since $\Q(u_p)=0$ for all $p \geq \bar p$ and $u_p \weak u_\infty$ in $W^{1,q}_0(\Om;\R^N)$ as $p_j\to\infty$, we deduce that $\Q(u_\infty)=0$, as desired. 

It remains to show that $u_{\infty}$ is indeed a minimiser of $\mathrm{F}_{\infty}$ in $\mathfrak{X}^{\infty}(\Om)$, and additionally that the energies converge. Fix an arbitrary $u\in \mathfrak{X}^{\infty}(\Om)$. By minimality and by noting that $\mathfrak{X}^{\infty}(\Om) \sub \mathfrak{X}^{p}(\Om)$ for any $p \in [\bar p ,\infty]$, we have the estimate
\[
\begin{split}
\mathrm{F}_{\infty}(u_{\infty})\, &= \, \smash{\lim_{q\to \infty}\mathrm{F}_{q}(u_{\infty})}\\
&\leq \, \liminf_{q\to \infty} \Big(\liminf_{p_j\to \infty}\mathrm{F}_{q}(u_{p})\Big)\\
&\leq \, \liminf_{p_j\to \infty}\mathrm{F}_{p}(u_{p})\\
&\leq \, \limsup_{p_j\to \infty}\mathrm{F}_{p}(u_p)\\
&\leq \, \limsup_{p_j\to \infty}\mathrm{F}_{p}(u)\\
& = \, \mathrm{F}_{\infty}(u),
\end{split}
\]
for any $u\in \mathfrak{X}^{\infty}(\Om)$. Hence $u_{\infty}$ is a minimiser of $\mathrm{F}_{\infty}$ over $\mathfrak{X}^{\infty}(\Om)$ and \eqref{1.4} follows. The choice $u=u_{\infty}$ in the above implies $\mathrm{F}_p(u_p)\larrow \mathrm{F}_{\infty}(u_{\infty})$, as $p_j\to \infty$. This completes the proof of the proposition and therefore of Theorem \ref{theorem1}.
\qed

\section{The equations for constrained minimisers in $L^p$ and  in $L^\infty$} \label{section3}

In this section we establish the proofs of Theorem \ref{theorem2} and Theorem \ref{theorem3}. The proof of the former is a relatively simple consequence of deep results in the generalised Kuhn-Tucker theory from \cite{Z}, whilst the proof of the latter is based on applying an appropriate liming process to the former result through compactness estimates.

\BPT \ref{theorem2}. Fix $p \in [\bar p, \infty)$. We begin with the simplifying observation that, the minimisation problem \eqref{1.14} can be rewritten as
\[
\frac{1}{p}\F_p(u_p)^p \, =\, \bigg\{ \frac{1}{p}\F_p(u)^p \ : \  u \in W^{1,\al p}_0(\Om;\R^N), \ \frac{1}{p}\G_p(u)^p- \frac{G^p}{p} \leq 0 \ \& \ \Q(u)=0 \bigg\}.
\]
This reformulation is a labour-saving device, drastically shortening the proof of this result. In view of assumption \eqref{1.15}, first we will show that the following functionals are Fr\'echet differentiable
\[
\begin{split}
\frac{1}{p}(\mathrm{F}_p)^p \ &: \ \  W^{1,\alpha p}_0\big(\Om;\R^N\big)\larrow \R,
\\
 \frac{1}{p}(\mathrm{G}_p)^p - \frac{G^p}{p} \ &: \ \  W^{1,\alpha p}_0\big(\Om;\R^N\big)\larrow \R.
\end{split}
\]
A direct computation gives the next formal expressions for their Gateaux derivatives
\[
\begin{split}
\bigg(  \d \Big[ \frac{1}{p}(\mathrm{F}_p)^p \Big]\bigg)_{\!\! u} (v) \, &=\, {\av}_{\!\!\!\Om}f[u]^{p-1}\big(f_{\eta}[u]\cdot v+ f_P[u]: \D v \big) \, \d {\mathcal{L}}^n,
\\
\bigg(  \d \Big[ \frac{1}{p}(\mathrm{G}_p)^p - \frac{G^p}{p}  \Big] \bigg)_{\!\! u}(v) \, &=\, {\av}_{\!\!\!\Om}g[u]^{p-1}\big(g_{\eta}[u]\cdot v+ g_P[u]: \D v \big)\, \d {\mathcal{L}}^n,
\end{split}
\]
for all $ u,v \in W^{1,\al p}_0(\Om;\R^N)$. We will now show the above formal expressions indeed define Fr\'echet derivatives, by employing relatively standard estimates through the H\"older inequality. We argue only for $\frac{1}{p}(\mathrm{F}_p)^p$, as the estimates for $\frac{1}{p}\big((\mathrm{G}_p)^p - {G^p}\big)$ are identical. Since $\al>1$ and $p\geq \bar p >n$, by Morrey's estimate we have
\[
\begin{split}
 \bigg|\, {\av}_{\!\!\!\Om}f[u] & ^{p-1}  \big(f_{\eta}[u]\cdot v+ f_P[u]: \D v \big) \, \d {\mathcal{L}}^n \bigg|  
\\
& \leq\,  {\av}_{\!\!\!\Om}|f[u]|^{p-1}\Big(|f_{\eta}[u]| | v| + |f_P[u]|  |\D v |\Big) \, \d {\mathcal{L}}^n
\\
& \leq\,  {\av}_{\!\!\!\Om}   C(\cdot, u)^p   \big(  1+  |\D u|^{\al}  \big)^{p-1}  \big(1+|\D u|^{\al-1}  \big) (|v|+|\D v|)\, \d {\mathcal{L}}^n
\\
& \leq\, 2^p \big\| |v| C(\cdot, u)^p\big\|_{L^\infty(\Om)} \,  {\av}_{\!\!\!\Om}   \Big(1+  |\D u  |^{\al  -1}  +  |\D u|^{\al p -\al}  + |\D u|^{\al p -1}  \Big) \, \d {\mathcal{L}}^n
\\
& \ \ + \, 2^p \big\| C(\cdot, u)^p\big\|_{L^\infty(\Om)}  \,  {\av}_{\!\!\!\Om}   \Big(1+  |\D u  |^{\al  -1}  +  |\D u|^{\al p -\al}  + |\D u|^{\al p -1}  \Big) |\D v|\, \d {\mathcal{L}}^n
\\
& \leq\, 2^p \big\| |v| C(\cdot, u)^p\big\|_{L^\infty(\Om)}   \Big(1+  \|\D u  \|_{L^{\al p}(\Om)}^{\al  -1}  +  \|\D u  \|_{L^{\al p}(\Om)}^{\al p -\al}  +  \|\D u  \|_{L^{\al p}(\Om)}^{\al p -1} \Big) 
\\
& \ \ + \, 2^p \big\| C(\cdot, u)^p\big\|_{L^\infty(\Om)}  \bigg(  \|\D v  \|_{L^{\al p}(\Om)} \\
&\  \ + \sum_{t \in \{ \al-1, \al p-\al, \al p-1 \}}   \|\D u  \|_{L^{ \frac{t\al p}{\al p -1}}(\Om)}^ t \|\D v  \|_{L^{\al p}(\Om)} \bigg) .
\end{split}
\]
We now show that the equations that the constrained minimiser satisfies take the form as given in \eqref{1.19} and \eqref{1.18}. Given the Fr\'echet derivatives and our assumption \eqref{1.16} on the range of $\d \Q$, we can invoke the generalised Kuhn-Tucker theory. By applying \cite[Theorem 48.B, p.p.\ 416-417]{Z} with (in the book's notation)
\[
\begin{split}
F_0  &:= \frac{1}{p}(\mathrm{F}_p)^p, \ \ \ F_1:= \frac{1}{p}(\mathrm{G}_p)^p - \frac{G^p}{p}, \ \ \ F_{3} := \Q,
\\
X &=U=N_2:= W^{1,\alpha p}_0\big(\Om;\R^N\big), \ \ \ Y:= \E, \ \ \ n=2,
\end{split}
\]
and by noting that $N_2$ herein is the entire (vector) space $W^{1,\alpha p}_0\big(\Om;\R^N\big)$,  we readily infer the claims made in \eqref{1.17}-\eqref{1.18}. The result ensues.
\qed
 \ms

Now we establish our last main result. 

\BPT \ref{theorem3}. The proof is divided into several steps.

\noi {\bf Step 1.} We first confirm that the measures defined by \eqref{1.20} are indeed finite, and show that their total variations are bounded uniformly in $p \in (\bar p, \infty)$. This will imply the convergence modes of \eqref{1.34} for some limiting $\mu_\infty, \nu_\infty \in \mM(\overline{\Om})$ along a subsequence $(p_j)_1^\infty$ as $j\to\infty$, as a consequence of the sequential weak* precompactness of bounded sets in the space of Radon measures. Indeed, if $\F_p(u_p)>0$, then since $f\geq0$ we have
\[
\begin{split}
\|\si_p\|(\overline{\Om}) \, &=\, \si_p(\overline{\Om})
\\
&=\, \frac{1}{\mL^n(\Om)}\int_{\Om} \frac{f(\cdot, \D u_p)^{p-1}}{\F_p(u_p)^{p-1}}\,\mathrm d \mL^n
\\
&=\, \frac{1}{\F_p(u_p)^{p-1}} \ {\av_{\Om}} f(\cdot, \D u_p)^{p-1}\,\mathrm d \mL^n
\\
&\leq\, \frac{1}{\F_p(u_p)^{p-1}} \bigg(\ {\av_{\Om}} f(\cdot, \D u_p)^{p}\,\mathrm d \mL^n\bigg)^{\!\! \frac{p-1}{p}}
\\
&=\, 1, \phantom{\Big|}
\end{split}
\]
whilst if $\F_p(u_p)=0$, then trivially $\|\si_p\|(\overline{\Om})=0$. In both cases, $\|\si_p\|(\overline{\Om}) \leq 1 $ for all $p \in (\bar p, \infty)$. The estimate for $\|\tau_p\|(\overline{\Om})$ is completely analogous, yielding $\|\tau_p\|(\overline{\Om}) \leq 1 $ for all $p \in (\bar p, \infty)$.

\ms

\noi {\bf Step 2.} By using assumption \eqref{1.27} and definition \eqref{1.20}, we have the following differential identity: for any fixed $v \in C^1_0\big(\overline{\Om};\R^N \big)$ and any $p\in (\bar p,\infty)$ we have
\[
\begin{split}
\int_\Om f(\cdot, \D v-\D u_p) \, \mathrm d \si_p  \, &=\, \int_\Om f(\cdot, \D v) \, \mathrm d \si_p - \int_\Om f(\cdot, \D u_p) \, \mathrm d \si_p 
\\
& \ \ + \, \int_\Om f_P(\cdot, \D u_p):(\D u_p - \D v) \, \mathrm d \si_p. 
\end{split}
\]
Indeed, by using that $f_P(x, P)\, =\, \A(x) \!:\! \big((\cdot)\ot P + P \ot (\cdot)\big)$,
we may compute
\[
\begin{split}
\int_\Om f(\cdot, \D v & - \D u_p)  \,  \mathrm d \si_p \, =\, \int_\Om \A:(\D v-\D u_p)\ot (\D v-\D u_p) \, \mathrm d \si_p
\\
&=\, \int_\Om \A: \D v \ot \D v \, \mathrm d \si_p \,-\,\int_\Om \A:\D u_p \ot \D u_p \, \mathrm d \si_p
\\
&\ \ + \int_\Om \A: \Big((\D u_p-\D v)\ot \D u_p + \D u_p \ot (\D u_p-\D v)\Big)\, \mathrm d \si_p
\\
&= \int_\Om f(\cdot, \D v) \, \mathrm d \si_p - \int_\Om f(\cdot, \D u_p) \, \mathrm d \si_p 
\, + \, \int_\Om f_P(\cdot, \D u_p):(\D u_p - \D v) \, \mathrm d \si_p. 
\end{split}
\]
We also note that the above established identity holds true over $\overline{\Om}$ as well, because $\si_p(\p\Om)=\tau_p(\p\Om)=0$.

\ms

\noi {\bf Step 3.} For any fixed $p\in(\bar p,\infty)$, by using \eqref{1.20}-\eqref{1.24} and \eqref{1.27}, we may rewrite \eqref{1.19} (obtained in Theorem \ref{theorem2}) as
\[
\begin{split}
 \La_p \int_{\overline{\Om}}  f_P(\cdot, \D u_p) : \D \phi  \, \d \si_p \, +\, \mathrm M_p \int_{\overline{\Om}} g_{\eta}(\cdot, u_p) \cdot \phi \, \d \tau_p
  \, =\, \big\langle \Psi_p , (\d \Q)_{u_p}(\phi) \big\rangle,
\end{split}
\]
for all test maps $\phi\in W^{1,\al p}_0(\Om;\R^N)$, whilst we also have that 
\[
\text{$\La_p \in [0,1]$, \ $\mathrm M_p \in [0,1]$ \ and \ $\Psi_p \in \bar \mB^{E^*}_1\!(0)$.}
\]
Further, by assumption \eqref{1.25}, the weak* topology of the dual space $\E^*$ is sequentially (pre-) compact on bounded sets. Thus, the previous memberships imply that, upon passing to a further subsequence as $j\to\infty$, symbolised again by $(p_j)_1^\infty$, there exist
\[
\text{$\La_\infty \in [0,1]$, \ $\mathrm M_\infty \in [0,1]$ \ and \ $\Psi_\infty \in \bar \mB^{E^*}_1\!(0)$,}
\]
such that the modes of convergence \eqref{1.33} hold true as $p_j\to\infty$. 

\ms

\noi {\bf Step 4.} By Steps 2 and 3, for $\phi := u_p -v$, where $v \in C^1_0\big(\overline{\Om};\R^N \big)$ is an arbitrary fixed map, for any fixed $p\in(\bar p,\infty)$ we have the identity
\[
\begin{split}
 \La_p   \int_{\overline{\Om}} & f(\cdot, \D u_p  - \D v) \,  \mathrm d \si_p    =\, \big\langle \Psi_p , (\d \Q)_{u_p}(u_p-v) \big\rangle \, 
 \\
 & -\, \mathrm M_p \int_{\Om} g_{\eta}(\cdot, u_p) \cdot (u_p -v) \, \d \tau_p \,+
 \, \La_p\bigg( \int_\Om f(\cdot, \D v) \, \mathrm d \si_p - \int_\Om f(\cdot, \D u_p) \, \mathrm d \si_p \bigg).
\end{split}
\]

\noi {\bf Step 5.} For any fixed $p\in(\bar p,\infty)$ and $v \in C^1_0\big(\overline{\Om};\R^N \big)$, we have the relations
\[
\left\{ \ \ 
\begin{split}
\int_{\overline{\Om}} f(\cdot, \D u_p-\D v) \, \mathrm d \si_p \, &\geq\, \al_0 \int_{\overline{\Om}} \big|\D u_p-\D v\big|^2 \, \mathrm d \si_p , 
\\
\int_{\overline{\Om}} f(\cdot, \D u_p) \, \mathrm d \si_p \, &=\, \F_p(u_p),
\end{split}
\right.
\]
where we have symbolised
\[
\al_0\,:=\, \min_{x\in \overline{\Om}} \bigg\{\min_{|Q|=1} \A(x): Q\ot Q \bigg\}\,>\,0.
\]
Let us first establish the claimed equality, beginning with the case that $\F_p(u_p)>0$. By definition \eqref{1.20} and assumption \eqref{1.27}, we may compute
\[
\begin{split}
\int_\Om f(\cdot, \D u_p) \, \mathrm d \si_p \, &=\, \frac{1}{\mL^n(\Om)}\int_{\Om} f(\cdot, \D u_p)\frac{f(\cdot, \D u_p)^{p-1}}{\F_p(u_p)^{p-1}}\,\mathrm d \mL^n
\\
&=\,  \frac{1}{\F_p(u_p)^{p-1}} \, {\av_{\Om}} f(\cdot, \D u_p)^p \,\mathrm d \mL^n
\\
&=\, \F_p(u_p), \phantom{\Big|}
\end{split}
\]
whilst for $\F_p(u_p)=0$ the equality follows trivially. To establish the claimed inequality, it suffices to note that by assumption \eqref{1.27} and by the variational representation of the minimum eigenvalue of the symmetric linear operator $\A(x) : \R^{N\by n} \larrow \R^{N\by n}$, we have that $\al_0>0$ and the inequality
\[
\al_0 |Q|^2\, \leq \, \A(x) : Q \ot Q
\]
for all $x\in \overline{\Om}$ and all $Q \in \R^{N\by n}$, where $|\cdot|$ is the Euclidean norm on $\R^{N\by n}$.

\ms

\noi {\bf Step 6.} By Steps 1, 3 and 5, and by using that $\F_p(u_p)\larrow \F_\infty(u_
\infty)$ as $p_j\to\infty$ (as shown in Theorem \ref{theorem1}), we may invoke Hutchinson's theory of measure-function pairs, in particular \cite[Sec.\ 4, Def.\ 4.1.1, 4.1.2, 4.2.1 and Th.\ 4.4.2]{H}, to infer that there exists a map
\[
V_\infty \ \in \, L^2\big(\overline{\Om},\si_\infty;\R^{N\by n}\big)
\]
such that, along perhaps a further subsequence $(p_j)_1^\infty$ we have
\[
\D u_p \si_p \, \weakstar \, V_\infty \si_\infty \ \text{ in }\mM\big(\overline{\Om};\R^{N\by n}\big),
\]
as $p_j\to \infty$, with the property that
\[
\int_{\overline{\Om}} \Phi(\cdot, V_\infty ) \, \d \si_\infty \, \leq \, \liminf_{p_j\to\infty} \int_{\overline{\Om}} \Phi(\cdot, \D u_p ) \, \d \si_p 
\]
for any $\Phi \in C\big( \overline{\Om}\by \R^{N\by n}\big)$ such that $\Phi(x,\cdot)$ is convex and of quadratic growth, for all $x\in \overline{\Om}$. Further, in view of assumptions \eqref{1.15}, \eqref{1.16}, \eqref{1.26}, \eqref{1.27}, and the modes of convergence established in Theorem \ref{theorem1} together with the convergence $\Psi_p \weakstar \Psi_\infty$ in $\E^*$
as $p_j\to \infty$, we have that
\[
\begin{split}
\big\langle \Psi_p, (\d \Q)_{u_p}(\phi) \big) \rangle & \larrow \big\langle \Psi_\infty , (\d \Q)_{u_\infty}(\phi) \big\rangle \ \ \ \,\text{ in }\R,
\\
\big(f_P(\cdot,\D u_p) :\D \phi \big) \si_p  & \,  \weakstar \, \big(f_P(\cdot,V_\infty):\D \phi\big) \si_\infty \ \text{ in }\mM\big(\overline{\Om}\big),
\end{split}
\]
for any fixed $\phi\in C^1_0\big(\overline{\Om};\R^N \big) \sub W^{1,\bar p}_0(\Om;\R^N)$. Hence, we may let $p_j\to \infty$ in Step 3, to deduce the equation
\[
\begin{split}
 \La_\infty \int_{\overline{\Om}}  f_P(\cdot, V_\infty) : \D \phi  \, \d \si_\infty \, +\, \mathrm M_\infty \int_{\overline{\Om}} g_{\eta}(\cdot, u_\infty) \cdot \phi \, \d \tau_\infty
  \, =\, \big\langle \Psi_\infty , (\d \Q)_{u_\infty}(\phi) \big\rangle,
\end{split}
\]
for any fixed $\phi\in C^1_0\big(\overline{\Om};\R^N \big)$. Further, by letting $p_j\to\infty$ in \eqref{1.18} we deduce that
\[
\mathrm M_\infty \big(\G_\infty(u_\infty) - G\big)=0.
\]
\noi {\bf Step 7.} The equations established in Step 6 will complete the proof of the theorem, upon establishing that 
\[
\ \ \ \ V_\infty \, =\, \D u^\star_\infty \ \ \ \si_\infty \text{-a.e.\ on }\overline{\Om},
\]
where $\D u^\star_\infty : \overline{\Om} \larrow \R^{N\by n}$ is some Borel measurable mapping which is a version of  $\D u_\infty \in L^\infty(\Om;\R^{N\by n})$, namely such that
\[
\ \ \ \ \D u_\infty \, =\, \D u^\star_\infty \ \ \ \mL^n\text{-a.e.\ on }\overline{\Om},
\]
(recall that $\p\Om$ is a nullset for the Lebesgue measure $\mL^n$). The remaining steps are devoted to establishing this claim, together with the approximability properties claimed in \eqref{1.32} for some sequence of mappings $(v_j)_1^\infty \sub C^1_0\big(\overline{\Om};\R^N \big)$, which will be constructed explicitly.
\ms

\noi {\bf Step 8.} If $\La_\infty=0$, then Step 6 completes the proof of Theorem \ref{theorem3} as the first term involving $V_\infty$ vanishes. Hence, we may henceforth assume that $\La_\infty >0$. Therefore, by passing perhaps to a further subsequence if necessary, we may assume that
\[
\ \ \ \ \ \La_{p_j} \, \geq \, \frac{\La_\infty}{2} \,> \, 0, \text{ for all }j\in\N.
\] 

\noi {\bf Step 9.} By Steps 3, 4, 5 and 8, the absolute continuity $\tau_p << \mL^n\LL_\Om$ and the bounds $0\leq \mathrm M_p \leq 1$ and $0\leq \tau_p\big(\overline{\Om}\big) \leq 1$, we have the estimate
\[
\begin{split}
\frac{\al_0\La_\infty}{2}  \int_{\overline{\Om}} \big|\D u_p  - \D v \big|^2 \,  \mathrm d \si_p \ \leq & \ \, \Big\langle \Psi_p , \, (\d \Q)_{u_p}(u_p)  -  (\d \Q)_{u_p}(v) \Big\rangle 
 \\
 & + \,\|u_p -v\|_{L^\infty(\Om)}\Big\{ \sup_{j\in\N} \, \big\|g_{\eta}(\cdot, u_{p_j})\big\|_{L^\infty(\Om)} \Big\}
 \\
& + \, \La_p\bigg(\int_{\overline{\Om}} f(\cdot, \D v) \, \mathrm d \si_p  \, - \, \F_p(u_p) \bigg),
\end{split}
\]
along the sequence $(p_j)_1^\infty$. By letting $j\to\infty$ in the above estimate, in view of Steps 1, 3 and 6 (for the choice $\Phi(x,Q):=|\D v(x)-Q|^2$) and assumption \eqref{1.26}, we infer that
\[
\begin{split}
\frac{\al_0\La_\infty}{2}  \int_{\overline{\Om}} \big|V_\infty  - \D v \big|^2 \,  \mathrm d \si_\infty \ \leq & \ \, \Big\langle \Psi_\infty  , \, (\d \Q)_{u_\infty}(u_\infty)  -  (\d \Q)_{u_\infty}(v) \Big\rangle  \\
 & + \,\|u_\infty -v\|_{L^\infty(\Om)}\Big\{ \sup_{j\in\N} \, \big\|g_{\eta}(\cdot, u_{p_j})\big\|_{L^\infty(\Om)} \Big\}
 \\
& + \, \La_\infty \bigg(\int_{\overline{\Om}} f(\cdot, \D v) \, \mathrm d \si_\infty  \, - \, \F_\infty(u_\infty)\bigg),
\end{split}
\]
for any fixed mapping $v\in C^1_0\big(\overline{\Om};\R^N \big)$.

\ms

\noi {\bf Step 10.} Let $(v_j)_1^\infty \sub C^1_0\big(\overline{\Om};\R^N \big)$ be any sequence of mappings satisfying the assumptions in \eqref{1.32}. We claim that there exists a subsequence of indices $(j_k)_1^\infty$ such that
\[
\left\{ \ \
\begin{array}{lll}
\D v_j \larrow V_\infty, & \text{ in }L^2\big(\overline{\Om},\si_\infty;\R^{N\by n}\big), &\text{ and $\si_\infty$-a.e.\ on }\overline{\Om}, 
\ms
\\
\D v_j \larrow \D u_\infty, & \text{ in }L^q\big({\Om};\R^{N\by n}\big), \, q\in [1,\infty), &\text{ and $\mL^n$-a.e.\ on }{\Om},
\end{array}
\right.
\]
as $j_k\to \infty$. Therefore, if $\D u^\star_\infty$ is defined as in \eqref{1.32}, then from the above we infer
\[
\left\{
\begin{split}
\D u^\star_\infty \, &=\, V_\infty, \ \ \ \, \si_\infty \text{-a.e.\ on }\overline{\Om},
\\
\D u^\star_\infty \, &=\, \D u_\infty, \ \ \mL^n\text{-a.e.\ on }{\Om},
\end{split}
\right.
\]
which completes the proof (subject to showing that at least one sequence of mapping $(v_j)_1^\infty$ with the desired properties exists). Let us now establish the above claims. If $(v_j)_1^\infty$ satisfies \eqref{1.32}, then by Step 5 and the resulting bound
\[
\al_0 \|\D v_j\|_{L^\infty(\Om)}^2 \, \leq\, \F_\infty (v_j) \, \leq\, \F_\infty(u_\infty) \,+ \, o(1)_{j\to\infty}
\]
in conjunction with the Vitaly convergence theorem, it follows that
\[
\ \ \ \ \ \ \left\{ \ \
\begin{array}{ll}
v_j \larrow u_\infty, & \text{ in }L^\infty ({\Om};\R^{N} ),
\ms
\\
\D v_j \larrow \D u_\infty, & \text{ in }L^q({\Om};\R^{N\by n}),\ q\in [1,\infty),
\end{array}
\right.
\]
along a subsequence of indices $(j_k)_1^\infty$ as $k\to\infty$. Consequently, by the estimate of Step 9, we have
\[
\begin{split}
\frac{\al_0\La_\infty}{2}  \int_{\overline{\Om}} \big|V_\infty  - \D v_j \big|^2 \,  \mathrm d \si_\infty \ \leq  & \ \, \Big\langle \Psi_\infty , \, (\d \Q)_{u_\infty}(u_\infty)  -  (\d \Q)_{u_\infty}(v_j) \Big\rangle 
 \\
 & + \,\|u_\infty -v_j\|_{L^\infty(\Om)}\Big\{ \sup_{j\in\N} \, \big\|g_{\eta}(\cdot, u_{p_j})\big\|_{L^\infty(\Om)} \Big\}
 \\
& + \, \La_\infty \bigg( \Big\{\sup_{\overline{\Om}} f(\cdot, \D v_j) \Big\}\,\si_\infty\big(\overline{\Om}\big)  \, - \, \F_\infty(u_\infty)\bigg).
\end{split}
\]
Since by Step 1 we have $0\leq \si_\infty\big(\overline{\Om}\big) \leq 1$, by noting that $\sup_{\overline{\Om}} f(\cdot, \D v_j) = \F_\infty (v_j)$ (due to the continuity of $\D v_j$ on $\overline{\Om}$) and that $v_j \larrow u_\infty$ strongly in $W^{1,\bar p}_0(\Om;\R^N)$, from the last estimate and assumption \eqref{1.32} we deduce that
\[
\begin{split}
\limsup_{j_k\to\infty}  \int_{\overline{\Om}} \big|V_\infty  - \D v_j \big|^2 \,  \mathrm d \si_\infty \, \leq  \, 0.
\end{split}
\]

\ms

\noi {\bf Step 11.} To complete the proof of Theorem \ref{theorem3}, it remains to show that at least one sequence of mapping $(v_j)_1^\infty \sub C^1_0\big(\overline{\Om};\R^N \big)$ exists, which satisfies the modes of convergence required by \eqref{1.32}. To this end we utilise (for the first time) the assumption that the bounded domain $\Om$ has Lipschitz boundary $\p\Om$, and we invoke the regularisation scheme introduced and utilised in the recent paper \cite{K4}. This method is based on results on the geometry of Lipschitz domains proved in Hofmann-Mitrea-Taylor \cite{HMT}, and is inspired by the regularisation schemes employed in Ern-Guermond \cite{EG}. If $\mathrm n\in L^\infty(\p\Om,\mH^{n-1};\R^n)$ be the outer unit normal vector field on $\p\Om$, by \cite[Sec.\ 2, 4]{HMT}, there exists a smooth vector field $\xi \in C^\infty_c(\R^n ;\R^n)$ which is globally transversal to $\mathrm n$ on $\p\Om$ with respect to the surface measure, namely exists $c>0$ such that
\[
\ \ \ \ \ \xi \cdot \mathrm  n  \geq c,  \quad \mH^{n-1}\text{-a.e.\ on }\p\Om.
\]
Additionally, $\xi$ can be chosen to satisfy $|\xi|\equiv 1$ in an open tubular neighbourhood $\{\dist(\cdot,\p\Om)<r\}$ around $\p\Om$ for some $r>0$, whilst vanishing on $\{\dist(\cdot,\p\Om)>2r\}$. (In the special case that $\p\Om$ happens to be a compact $C^\infty$ manifold, then we can simply choose $\xi$ to be a smooth extension of $\mathrm n$, and the transversality condition is satisfied with $c=1$.) Further, it can be shown that there exists $\e_0, \ell>0$ such that, for all $\e \in (0,\e_0)$ we have
\[
\ \ \ \inf_{x\in\p\Om} \dist\big(x+\e \ell\xi(x),\, \p\Om \big) \geq \, 2\e.
\]
We now define our adapted global mollifiers, taken from \cite{K4}. Let us select any function $\varrho \in C^\infty_c(\mB_1(0))$ which satisfies $\varrho\geq 0$ and $\|\varrho\|_{L^1(\R^n)}=1$ (for instance the ``standard" mollifying kernel as in \cite{KV}). For any $v\in L^\infty(\Om;\R^N)$, extended to $\R^n\set \Om$ by zero, we define for any $\e \in (0,\e_0)$ the map $\K^\e v : \R^n \larrow \R^N$, by setting
\[
(\K^\e v)(x) \,:=\, \int_{\R^n}  v\big(x+\e \ell\xi(x)-\e y \big)\, \varrho(y)\, \d y.
\]
Then, by \cite[Prop.\ 12, p. 18]{K4}, for any $u \in W^{1,\infty}_0(\Om;\R^N)$ and $\e \in (0,\e_0)$:

\ms

\noi $\bullet$ We have that $\K^\e u \in C^\infty_0\big( \Om ; \R^N \big)$, and the identity
\[
\D(\K^\e u) \, =\, \K^\e (\D u)   \,+\, \e\ell  (\K^\e (\D u) )(\D\xi)^\top,
\]
everywhere on $\overline{\Om}$.
\ms

\noi $\bullet$ We have $\K^\e u  \larrow u$ in $W^{1,q}_0(\Om;\R^N)$ for all $q\in[1,\infty)$, and in $C^\ga(\overline{\Om};\R^N)$ for all $\ga \in (0,1)$, as $\e\to 0$. Further, $\K^\e u \, \weakstar \, u $ in $W^{1,\infty}_0(\Om;\R^N)$, as $\e\to 0$.
\ms

\noi $\bullet$  For any $\Theta \in C\big( \overline{\Om} \by \R^{N\by n}\big)$, satisfying for any $x\in \overline{\Om}$ that $\Theta(x,\cdot)$ is convex on $\R^{N\by n}$ with $\Theta(x,\cdot) \geq \Theta(x,0)=0$, and also that the partial derivative $\Theta_P$ exists and is continuous on $\overline{\Om} \by \R^{N\by n}$, we will show that there exists a modulus of continuity $\om \in C\big([0,\infty);[0,\infty)\big)$ with $\om(0)=0$ which is independent of $\e$ and $x$, such that 
\[
\Theta\big(x, \D(\K^\e u )(x)\big) \,\leq\, \underset{\Om\cap \mB_\e(x+\e \ell\xi(x))}{\ess\sup}\, \Theta(\cdot, \D u ) \,+\,\om(\e),
\]
for any $x\in \Om$. In particular, the above estimate implies that
\[
\big\|\Theta\big(\cdot, \D(\K^\e u )\big) \big\|_{L^\infty(\Om)} \, \leq  \, \|\Theta(\cdot, \D u )\|_{L^\infty(\Om)} \,+\,\om(\e).
\]
(This was already established in \cite{K4}, but without $x$-dependence for $\Theta$.) As a result, 
\[
\limsup_{\e\to0} \big\|\Theta\big(\cdot, \D(\K^\e u )\big) \big\|_{L^\infty(\Om)} \, \leq  \, \|\Theta(\cdot, \D u )\|_{L^\infty(\Om)} .
\]
Let us now establish the claimed estimate. For any fixed $R>0$ such that 
\[
R \, > \, \big(\|\D\xi \|_{L^\infty(\R^n)}+1\big) \big(\|\D u \|_{L^\infty(\Om)}+1\big), 
\]
we set
\[
\om(t)\,:=\, t \ell R \|\Theta_P\|_{C(\overline{\Om}\by \bar \mB_R(0))} \, + \underset{|z-x|\leq t \ell R}{\sup_{|Q|\leq R}} \Big| \Theta\big(x,Q\big) - \Theta\big(z,Q\big) \Big|, \ \ \ t\geq 0.
\]
Then, we have
\[
\begin{split}
\Theta \big(x,\D(\K^\e u) (x)\big)  & = \Theta \Big(x,\K^\e (\D u)(x)  \,+\, \e\ell \K^\e (\D u)(x)(\D\xi(x))^\top \Big) 
\\
& \leq  \Theta \big(x,\K^\e (\D u)(x)\big) + \|\Theta_P(x,\cdot) \|_{C(\bar \mB_R(0))} \big\|\e\ell  \K^\e (\D u)(\D\xi)^\top \big\|_{L^\infty(\Om)}
\\
& \leq  \Theta \big(x,\K^\e (\D u)(x)\big) + \e\ell \|\Theta_P\|_{C(\overline{\Om}\by \bar \mB_R(0))} \|\D u \|_{L^\infty(\Om)} \|\D\xi \|_{L^\infty(\R^n)} 
\\
&\leq  \Theta \big(x,\K^\e (\D u)(x)\big)\, +\, \e \ell R \|\Theta_P\|_{C(\overline{\Om}\by \bar \mB_R(0))},
\end{split}
\]
for any $x\in\Om$. Further, since $\varrho \mL^n$ is a probability measure on $\R^n$, by Jensen's inequality, we have
\[
\begin{split}
 \Theta \big(x,\K^\e (\D u)(x)\big) \, & =\, \Theta \left(x,\int_{\R^n}  \D u\big(x+\e \ell\xi(x)-\e y \big)\, \varrho(y)\, \d y \right)
\\
& \leq\, \int_{\R^n} \Theta \Big(x, \D u \big(x+\e \ell\xi(x)-\e y \big)\Big) \varrho(y)\, \d y
\\
&\leq\, \underset{y\in \mB_1(0)}{\ess\sup}\, \Theta \Big(x,\D u\big(x+\e \ell\xi(x)-\e y \big)\Big) 
\\
&=\, \underset{z\in \mB_\e(x+\e \ell\xi(x))}{\ess\sup}\, \Theta\big(x,\D u(z)\big) 
\\
&\leq \, \underset{z\in \mB_\e(x+\e \ell\xi(x))}{\ess\sup}\, \Theta\big(z,\D u(z)\big) \\
&\ \ \ + \underset{|z-x|\leq \e \ell \|\D \xi\|_{L^\infty(\R^n)}}{\sup_{|Q|\leq \|\D u\|_{L^\infty(\Om)}}} \Big| \Theta\big(x,Q\big) - \Theta\big(z,Q\big) \Big|,
\end{split}
\]
for any $x\in\Om$. By using that $\D u \equiv 0$ on $\R^n \set \Om$ and our assumptions on $\Theta$, the previous two estimates yield that
\[
\begin{split}
\Theta \big(x,\D(\K^\e u) (x)\big)  \, & \leq\, \underset{\mB_\e(x+\e \ell\xi(x))}{\ess\sup}\, \Theta (\cdot,\D u )\, +\, \om(\e)
\\
& =\, \max\bigg\{ \underset{\Om\cap \mB_\e(x+\e \ell\xi(x))}{\ess\sup}\,\Theta (\cdot,\D u ),\ \underset{\mB_\e(x+\e \ell\xi(x)) \set\Om}{\ess\sup}\, \Theta (\cdot,\D u ) \bigg\}\, +\, \om(\e)
\\
&=\, \max\bigg\{ \underset{\Om\cap \mB_\e(x+\e \ell\xi(x))}{\ess\sup}\, \Theta (\cdot,\D u ),\ \underset{\mB_\e(x+\e \ell\xi(x)) \set\Om}{\ess\sup}\, \Theta (\cdot,0) \bigg\}\, +\, \om(\e)
\\
&=\, \underset{\Om\cap \mB_\e(x+\e \ell\xi(x))}{\ess\sup}\, \Theta (\cdot,\D u )\, +\, \om(\e),
\end{split}
\]
for any $x\in\Om$. 

In conclusion, the above establish the existence of an approximating sequence as in \eqref{1.32}, by taking 
\[
\ \ \ \ v_j \, := \, \K^{\e_j} u_\infty  \in \, C^1_0\big( \overline{\Om} ; \R^N \big),\ \ j \in \N,
\]
along any infinitesimal sequence $(\e_j)_1^\infty$ satisfying $\e_j\to 0$, for the choice $\Theta := f$, which is admissible because of our hypotheses \eqref{1.5}, \eqref{1.6} and \eqref{1.27}.

\ms

\noi {\bf Step 12.} By putting together all the previous, the proof of Theorem \ref{theorem3} has now been completed.
\qed
\ms

\section{Explicit classes of nonlinear operators}  \label{section5}

In this section we provide various examples of nonlinear operators $\Q$ as in \eqref{1.2}, satisfying our assumptions. Note that each our main results, Theorems \ref{theorem1}, \ref{theorem2} and \ref{theorem3}, have been obtained with progressively stronger assumption on the operator $\Q$ which expresses one of constraints in the admissible class. For the sake of clarity, in the next table we list in a concise way which assumptions are required to be satisfied by $\Q$, in order to obtain the corresponding result (assuming that $f, g$ satisfy separately their respective required assumptions, the table concerns $\Q$ solely).
\[
\underset{\text{Table 1. The list of assumptions on $\Q$ to obtain the main results.}}{\left|
\!\begin{array}{rrrr}
\hline
\phantom{\Big|}\eqref{1.9} &\phantom{a} &  \phantom{a}  & \Mapsto \text{Theorem \ref{theorem1}}
 \\
\hline
\phantom{\Big|} \eqref{1.9} & \& \ \ \eqref{1.16} & \phantom{a} & \Mapsto \text{Theorem \ref{theorem2}}
 \\
 \hline
\phantom{\Big|} \eqref{1.9} & \& \ \ \eqref{1.16}  &  \& \ \ \eqref{1.25}\ \ \& \ \ \eqref{1.26} & \Mapsto  \text{Theorem \ref{theorem3}}
 \\
\hline
\end{array}
\!
\right|}
\]

\subsection{Pointwise constraints, unilateral constraints and inclusions} \label{S4.1} The nonlinear operator of \eqref{1.2} we are using in the admissible class of \eqref{1.4}, can include the following model cases:

\ms

\noi \textbf{Case 1.} $\Pi(x,u(x))=0$ for a.e.\ $x\in \Om$, where $\Pi : \Om \by \R^N \larrow \R^M$ is given.
\ms

\noi \textbf{Case 2.} $\Pi(x,u(x))\leq 0$ for a.e.\ $x\in \Om$, where  $\Pi : \Om \by \R^N \larrow \R$  is given.
\ms

\noi \textbf{Case 3.} $u(x) \in \mathcal{K}$ for a.e.\ $x\in \Om$, where  $\mathcal{K} \sub \R^N$  is a given closed set.
\ms

\noi Constraints as in Case 1 are sometimes called \emph{holonomic} (see for instance \cite{AB}). We now elaborate on the assumptions required to be fulfilled in each of these cases.

\ms

\begin{proposition}[Case 1] \label{Prop7} Let $\Pi \in C^1(\overline{\Om} \by \R^N;\R^M)$, $M\in\N$. By defining  
\beq
\label{4.1}
\Q \ : \ \ W^{1,\bar p}_0 (\Om;\R^N) \larrow L^1(\Om;\R^M), \ \ \ \ \Q(u) \,:= \, \Pi(\cdot,u),
\eeq
and setting $\emph \E := L^1(\Om;\R^M)$, we have the following:
\ms

\noi \emph{(i)} The zero set of $\Q$ equals
\[
\Q^{-1}(\{0\}) \, =\, \Big\{ v \in W^{1,\bar p}_0 (\Om;\R^N) \ :\ \Pi(x,v(x)) = 0, \text{ a.e. }x\in \Om\Big\},
\]
and assumption \eqref{1.9} is always satisfied.

\ms

\noi \emph{(ii)}  If for any $x\in \Om$ we have
\[
\big\{ \Pi(x,\cdot )=0\big\} \, \sub  \, \big\{ \Pi_\eta (x,\cdot )=0\big\}, 
\]
namely when all points in the zero set are critical points, then $\Q$ satisfies \eqref{1.16}. 
 
\ms

\noi \emph{(iii)} Assumptions \eqref{1.25} and \eqref{1.26} are always satisfied. 
\end{proposition}

The choice of $\E$ is deliberately made ``as large as possible", as then the Lagrange multipliers of Theorems \ref{theorem2} and \ref{theorem3} are valued in the smaller space $\E^* = L^\infty(\Om;\R^M)$.

\BPP \ref{Prop7}. (i) Follows directly from the definitions, by the continuity of $\Pi$ and by Morrey's estimate, because $\bar p>n$. 
\smallskip

\noi (ii) Indeed, since
\[
(\mathrm d \Q)_u(\phi) \, =\, \Pi_\eta(\cdot, u) \cdot \phi,
\]
if $u \in \Q^{-1}(\{0\})$, then $\Pi(\cdot, u)=0$ a.e.\ on $\Om$ and therefore  $\Pi_\eta(\cdot, u)=0$ a.e.\ on $\Om$, which implies that $(\mathrm d \Q)_u =0$, hence its image is the closed trivial subspace $\{0\} \sub L^1(\Om;\R^M)$.   

\smallskip

\noi (iii) Note first that $L^1(\Om;\R^M)$ is separable. Also, if we have $u_m \weak u$  and $\phi_m \weak \phi$ in $W^{1,\bar p}_0 (\Om;\R^N)$ as $m\to\infty$, then by Morrey's theorem and the compactness of the imbedding of H\"older spaces we have $u_m\larrow u$ and also $\phi_m \larrow \phi$ in $C\big(\overline{\Om};\R^N\big)$ as $m\to\infty$. Hence, we have as $m\to\infty$ that
\[
(\mathrm d \Q)_{u_m}(\phi_m) =\Pi_\eta(\cdot, u_m)\cdot \phi_m \, \larrow \, \Pi_\eta(\cdot, u)\cdot \phi = (\mathrm d \Q)_u(\phi),
\]
in $C\big(\overline{\Om};\R^M\big)$, which a fortiori implies strong convergence in $L^1(\Om;\R^M)$.  \qed
\ms

We note that the proof of (iii) above is immediate if one assumes the additional hypothesis of (ii), since then $(\mathrm d \Q)_{u_m}=0$ for any sequence $(u_m)_1^\infty \sub \Q^{-1}(\{0\})$.

\begin{proposition}[Case 2] \label{Prop8} Let $\Pi \in C^1(\overline{\Om} \by \R^N)$ and let us define $\pi : \R \larrow \R$ by
\beq 
\label{4.2}
\pi(t)\,:=\,
\left\{
\begin{array}{ll}
0, & t \leq 0,
\\
t^2, & t>0.
\end{array}
\right.
\eeq
By defining the operator
\[
\Q \ : \ \ W^{1,\bar p}_0 (\Om;\R^N) \larrow L^1(\Om), \ \ \ \ \Q(u) \,:= \, \pi(\Pi(\cdot, u)),
\]
for $\emph \E := L^1(\Om)$, we have the following:

\ms

\noi \emph{(i)} The zero set of $\Q$ equals
\[
\Q^{-1}(\{0\}) \, =\, \Big\{ v \in W^{1,\bar p}_0 (\Om;\R^N) \ :\ \Pi(x, v(x)) \leq 0, \text{ a.e. }x\in \Om\Big\},
\]
and assumption \eqref{1.9} is always satisfied.

\ms

\noi \emph{(ii)}  If for any $x\in \Om$ it holds that
\[
\big\{ \Pi(x,\cdot )=0\big\} \, \sub  \, \big\{ \Pi_\eta (x,\cdot )=0\big\},
\]
then $\Q$ satisfies assumption \eqref{1.16}.
\ms

\noi \emph{(iii)} Assumptions \eqref{1.25} and \eqref{1.26} are always satisfied. 
\end{proposition}

\BPP \ref{Prop8}. (i) Follows as in the proof of Proposition \ref{Prop7}(i), upon noting that $\{\pi \leq 0\}=(-\infty,0]$.

\smallskip

\noi (ii) Since
\[
(\mathrm d \Q)_u(\phi) \, =\, \pi'\big(\Pi(\cdot, u)\big)\Pi_\eta(\cdot, u) \cdot \phi,
\]
if $u \in \Q^{-1}(\{0\})$, then $\Pi(\cdot, u)\leq 0$ a.e.\ on $\Om$ and therefore  $\pi'(\Pi(\cdot, u))=0$ a.e.\ on $\Om$ because $\{\pi'=0\}=(-\infty, 0]$, which implies that $(\mathrm d \Q)_u =0$, hence its image is the closed trivial subspace $\{0\} \sub L^1(\Om;\R^M)$ and \eqref{1.16} is satisfied.
\smallskip

\noi (iii) Similar to the proof of Proposition \ref{Prop7}(iii), using the different expression for the differential $\d \Q$ as above.
\qed
\ms

\begin{proposition}[Case 3] \label{Prop9}  Let $\mathcal{K} \sub \R^N$ be a closed set with $\mathcal{K} \neq \emptyset$. Then, there exists $\Pi \in C^\infty(\R^N)$ satisfying $\mathcal{K} = \{\Pi=0\} \sub \{\Pi_\eta =0\}$. Further, if one defines
\[
\Q \ : \ \ W^{1,\bar p}_0 (\Om;\R^N) \larrow L^1(\Om), \ \ \ \ \Q(u) \,:= \, \Pi(u),
\]
and $\emph \E := L^1(\Om)$, then we have 
\[
\Q^{-1}(\{0\}) \, =\, \Big\{ v \in W^{1,\bar p}_0 (\Om;\R^N) \ :\ v(x) \in \mathcal{K}, \text{ a.e. }x\in \Om\Big\},
\]
and $\Q$ satisfies \eqref{1.9}, \eqref{1.16}, \eqref{1.25} and \eqref{1.26}.
\end{proposition}

\BPP \ref{Prop9}. It is well-known that for every such set $\mathcal{K}$, there exists a function $\Pi \in C^\infty(\R^N)$ with the claimed properties. A proof of this fact can be found e.g.\ in \cite[Sec.\ 1.1.13, p.\ 25]{N} (the claimed inclusion is not explicitly stated, but follows from the method of proof by the smooth Urysohn lemma). The rest follows from Proposition \ref{Prop7}.
\qed
\ms


\subsection{Integral and isoperimetric constraints} \label{S4.2} The nonlinear operator of \eqref{1.2} can also cover the following important case of constraint:
\[
\int_\Om h(\cdot, u,\D u) \, \d \mL^n \leq H,
\]
when $h : \Om \by \R^N \by \R^{N\by n}\larrow \R$ and $H \in \R$ are given.

\begin{proposition} \label{Prop10} Let $h : \Om \by \R^N \by \R^{N\by n}\larrow \R$ satisfy the assumptions \eqref{1.5}-\eqref{1.7} and \eqref{1.15} that $f, g$ are assumed to satisfy, with $\al \leq \bar p$. Let also $H \in \R$ be given and let $\pi :\R \larrow \R$ be as in \eqref{4.2}. Then, by defining the operator
\[
\Q \ : \ \ W^{1,\bar p}_0 (\Om;\R^N) \larrow \R, \ \ \ \ \Q(u) \,:= \, \pi\left(\, \int_\Om h(\cdot, u,\D u)\, \d \mL^n -H \! \right),
\]
and setting $\emph \E := \R$, we have the following:

\ms

\noi \emph{(i)} The zero set of $\Q$ equals
\[
\Q^{-1}(\{0\}) \, =\, \left\{ v \in W^{1,\bar p}_0 (\Om;\R^N) \ :\ \int_\Om h(\cdot, v,\D v)\, \d \mL^n \leq H \right\}
\]
and assumption \eqref{1.9} is satisfied.

\ms

\noi \emph{(ii)} $\Q$ satisfies \eqref{1.16}, \eqref{1.25} and \eqref{1.26}. 
\end{proposition}

\BPP \ref{Prop10}. (i) If $\Q(u_m)=0$ and $u_m \weak u$ in $ W^{1,\bar p}_0 (\Om;\R^N)$ as $m\to\infty$, then since $\{\pi \leq 0\}=(-\infty, 0]$, we have
\[
\int_\Om h(\cdot, u_m,\D u_m) \, \d \mL^n  - H \, \leq \, 0.
\]
Since $h$ satisfies \eqref{1.5}-\eqref{1.7} for $\al\leq \bar p$, by standard results (see e.g.\ \cite{D}), the functional $u \mapsto \|h(\cdot, u,\D u) \|_{L^1(\Om)}$ is weakly lower-semicontinuous in $W^{1,\bar p}_0(\Om;\R^N)$. Hence
\[
\int_\Om h(\cdot, u,\D u) \, \d \mL^n  - H \, \leq\, \liminf_{m\to\infty} \int_\Om h(\cdot, u_m,\D u_m) \, \d \mL^n  - H \, \leq \, 0.
\]
Therefore, $\Q(u)=0$, yielding that $\Q^{-1}(\{0\})$ is weakly closed and hence \eqref{1.9} is satisfied.

\smallskip

\noi (ii) By a computation, the Gateaux derivative of $\Q$ is given by
\[
(\mathrm d \Q)_u(\phi) = \pi' \! \left(\int_\Om h(\cdot, u,\D u)\, \d \mL^n -H \! \right) \!\! \int_\Om \Big[ h_\eta(\cdot, u,\D u)\cdot \phi \,+\, h_P(\cdot, u,\D u): \D \phi \Big]\, \d \mL^n ,
\]
and assumption \eqref{1.15} for $h$ implies that $\d \Q$ is (jointly) continuous on $W^{1,\bar p}_0(\Om;\R^N) \by W^{1,\bar p}_0(\Om;\R^N)$. Further, if $u \in \Q^{-1}(\{0\})$, then by part (i) we have
\[
\int_\Om h(\cdot, u,\D u) \, \d \mL^n -H \, \leq \, 0,
\]
and therefore the first factor of $(\mathrm d \Q)_u(\phi)$ vanishes because $\{\pi'=0\}=(-\infty, 0]$. Thus, $(\mathrm d \Q)_u =0$ when $u \in \Q^{-1}(\{0\})$, and hence its image is the closed trivial subspace $\{0\} \sub \R$, yielding that \eqref{1.16} is satisfied.
\smallskip

\noi (iii) For any sequences $u_m \weak u$ in $\Q^{-1}(\{0\}) \sub W^{1,\bar p}_0 (\Om;\R^N)$ and $\phi_m \weak \phi$ in $W^{1,\bar p}_0 (\Om;\R^N)$ as $m\to\infty$, by part (ii) we have 
\[
(\mathrm d \Q)_{u_m}(\phi_m)\, =\, 0 \larrow 0 \, =\, (\mathrm d \Q)_{u}(\phi)
\]
as $m\to\infty$, hence \eqref{1.25} and \eqref{1.26} are satisfied.
 \qed


\subsection{Quasilinear second order differential constraints} \label{S4.3}

The operator $\Q$ of \eqref{1.2} can also cover the case of various types of nontrivial PDE constraints. As an example, we discuss the case of quasilinear divergence second order systems of PDE of the form
\beq
\label{4.3}
\div \big( A(\cdot,u,\D u) \big) \, =\, B(\cdot,u,\D u) \ \ \text{ in }\Om,
\eeq
where the coefficients maps $A : \Om \by \R^N \by \R^{N\by n} \larrow \R^{N\by n}$ and $B : \Om \by \R^N \by \R^{N\by n} \larrow \R^{N}$ are given. Given the plethora of possibilities on the assumptions for such systems, the discussion in this subsection is less formal and is only aimed as a general indication of the admissible choices for $\Q$.

Suppose that $A,B$ are $C^1$ and satisfy appropriate growth bounds, and also that $P\mapsto A(\cdot,\cdot,P)$ a monotone map, and that the set of weak solutions to the system \eqref{4.3} is strongly precompact in $W^{1,\bar p}_0 (\Om;\R^N)$. A sufficient conditions for strong precompactness in $W^{1,\bar p}_0 (\Om;\R^N)$ for the set of weak solutions is for example a global $C^{1,\ga}$ or a $W^{2,1+\ga}$ a priori uniform bound on the set of solutions, for some $\ga \in (0,1)$. Appropriate assumptions on the coefficients $A,B$ that allow the derivation of such a priori bounds can be found e.g.\ in \cite{GT} for $N=1$ and in \cite{GM} for $N\geq2$. Then, by defining the operator
\[
\Q \ : \ W^{1,\bar p}_0 (\Om;\R^N) \larrow W^{-1,\bar p'} (\Om;\R^N)
\]
as
\[
\langle \Q(u),\psi\rangle \,:=\, \int_\Om \Big[A(\cdot,u,\D u):\D \psi \,+ \, B(\cdot,u,\D u) \cdot \psi \Big]\, \d \mL^n,
\]
and setting also $\E := W^{-1,\bar p'} (\Om;\R^N)$, assumptions \eqref{1.9},  \eqref{1.16},  \eqref{1.25} and  \eqref{1.26} are satisfied, with
\[
\Q^{-1}(\{0\}) \, =\, \bigg\{u \in W^{1,\bar p}_0 (\Om;\R^N)\ : \ \div \big( A(\cdot,u,\D u) \big)  = B(\cdot,u,\D u) \text{ weakly in }\Om \bigg\}.
\]
Note first that the expression of $\Q^{-1}(\{0\})$ is immediate by the definition of the differential operator $\Q$. Next, note that by assumption, for any sequence of weak solutions $(u_m)_1^\infty \sub W^{1,\bar p}_0 (\Om;\R^N)$ to \eqref{4.3}, there exists $u \in W^{1,\bar p}_0 (\Om;\R^N)$ such that $u_m \larrow u$ strongly along a subsequence $m_j\to\infty$. By applying this to any sequence $(u_m)_1^\infty \sub \Q^{-1}(\{0\})$ (namely sequence of solutions) for which $u_m \weak u$ as $m\to \infty$, by passing to the limit in the weak formulation for fixed $\psi \in W^{1,\bar p}_0 (\Om;\R^N)$, which reads
\[
\int_\Om \Big[A(\cdot,u_m,\D u_m):\D \psi \,+ \, B(\cdot,u_m,\D u_m) \cdot \psi \Big]\, \d \mL^n \, =\, 0,
\]
we get that $u\in \Q^{-1}(\{0\})$, as the convergence is in fact strong. Hence, \eqref{1.9} is satisfied. Further, under appropriate bounds, the operator $\Q$ is Fr\'echet differentiable and
\[
\begin{split}
\big\langle (\d \Q)_u(\phi),\psi \big\rangle \, &=\, \int_\Om \Big[  A_\eta(\cdot,u,\D u) \cdot \phi \,+\, A_P(\cdot,u,\D u) : \D\phi  \Big] : \D \psi \, \d \mL^n
 \\
&+\, \int_\Om \Big[  B_\eta(\cdot,u ,\D u ) \cdot \phi \,+\, B_P(\cdot,u ,\D u ) : \D\phi  \Big] \cdot \psi \, \d \mL^n.
\end{split}
\]
To see that the image of $(\d \Q)_u : W^{1,\bar p}_0 (\Om;\R^N) \larrow \E$ is closed for any fixed $u\in \Q^{-1}(\{0\})$, let $(T_m)_1^\infty \sub \mathrm{Rg}\big( (\d \Q)_u\big) \sub \E$ be a sequence in the range with $T_m \larrow T$ strongly in $\E$ as $m\to \infty$. Since $T_m \in \mathrm{Rg}\big( (\d \Q)_u\big)$, exists $\phi_m \in W^{1,\bar p}_0 (\Om;\R^N)$ solving the following linear second order system
\[
\begin{split}
- \div \Big(  A_\eta(\cdot,u,\D u) \cdot \phi_m \,& +\, A_P(\cdot,u,\D u):\D\phi_m \Big)
 \\
 & +\, B_P(\cdot,u ,\D u ) : \D\phi_m\,+\, B_\eta(\cdot,u ,\D u ) \cdot \phi_m  \,= \, T_m.
\end{split}
\]
By the monotonicity of the above system (due to our earlier assumption), under appropriate conditions one has a uniform bound in $W^{1,\bar p}_0 (\Om;\R^N)$, yielding the weak compactness of the sequence of solutions $(\phi_m)_1^\infty$, which establishes the closedness of $\mathrm{Rg}\big( (\d \Q)_u\big) \sub \E$ and \eqref{1.16} ensues.

Finally, for any sequence $(u_m)_1^\infty\sub \Q^{-1}(\{0\})$ satisfying $u_m \weak u$ as $m\to \infty$ and any $\phi \in W^{1,\bar p}_0 (\Om;\R^N)$, there exists $m_j\to\infty$ such that $u_m \larrow u$ as $m_j\to \infty$. These facts imply that $(\d \Q)_{u_m}(u_m)\larrow (\d \Q)_u(u)$ and also $(\d \Q)_{u_m}(\phi)\larrow (\d \Q)_u(\phi)$, both strongly in $ \E$ as $m\to\infty$. Hence, \eqref{1.25} and \eqref{1.26} are satisfied.

\subsection{Null Lagrangians and determinant constraints} \label{S4.4}

We close this paper with the observation that Theorem \ref{theorem1} holds true even when $\Q$ expresses a fully nonlinear pointwise jacobian determinant constraint, or even a more general pointwise PDE constraint driven by a null Lagrangian. As an explicit example, let $n=N$ and consider the differential operator
\[
\Q \ : \ \ W^{1,\bar p}_0 (\Om;\R^n) \larrow W^{-1,(\bar p/n)'} (\Om),
\]
by setting
\[
\Q(u)\, :=\, \det(\D u) - h,
\]
for a fixed $h \in L^{\bar p/n} (\Om)$, satisfying the necessary compatibility condition 
\[
\int_\Om h\, \d \mL^n\,=\,0. 
\]
We also take 
\[
\E \, := \, W^{-1,(\bar p/n)'} (\Om) \, = \, \big( W_0^{1,\bar p/n}(\Om) \big)^*. 
\]
Then, we have
\[
\Q^{-1}(\{0\}) \, =\, \Big\{u \in W^{1,\bar p}_0 (\Om;\R^n)\ : \ \det(\D u) = h \text{ a.e.\ in }\Om \Big\}.
\]
It follows that \eqref{1.9} is satisfied by the well-known property of weak continuity for jacobian determinants (see e.g.\ \cite[Th.\ 8.20, p.\ 395]{D}). However, the situation is more complicated regarding the satisfaction of the remaining assumptions.
If additionally $n=2$, then \eqref{1.25} and \eqref{1.26} are also satisfied. Indeed, since
\[
(\d \Q)_u(\phi) \, =\, \mathrm{cof}(\D u) :\D \phi,
\]
and since for $u=\phi$ we have the identity
\[
(\d \Q)_u(u) \, =\, \mathrm{cof}(\D u) :\D u \, =\, n\det(\D u),
\]
for any $(u_m)_1^\infty \sub \Q^{-1}(\{0\})$ with $u_m \weak u$ as $m\to\infty$, we have
\[
\begin{split}
(\d \Q)_{u_m}(u_m) \, =\, n \det(\D u_m) \, \weak \, n \det(\D u) \, =\,  (\d \Q)_{u}(u) 
\end{split}
\]
in $L^{\bar p/2} (\Om)$ as $m\to\infty$, whilst for any $\phi \in W^{1,\bar p}_0 (\Om;\R^2)$ we have
\[
(\d \Q)_{u_m}(\phi) \weak (\d \Q)_{u}(\phi)
\]
in $L^{\bar p/2} (\Om)$ as $m\to\infty$, by the linearity of the cofactor operator when $n=2$. Then, the compactness of the imbedding
\[
L^{\bar p/2} (\Om) \, \Subset \, W^{-1,(\bar p/2)'} (\Om) 
\]
implies that the above modes of convergence are in fact strong in $\E = W^{-1,(\bar p/2)'} (\Om)$. However, it is not clear when assumption \eqref{1.16} is satisfied, or when \eqref{1.26} is satisfied in the case that $n\geq 3$. This means Theorems \ref{theorem2} and  \ref{theorem3} as they stand do not apply to the case of jacobian constraints. This does \emph{not} mean that it is impossible to derive the associated equations, it merely means that in this case of such a highly nonlinear constraint a different specialised method of proof is required.


\ms

\bibliographystyle{amsplain}

\end{document}